\documentclass[10pt]{amsart}
\usepackage{latexsym}
\usepackage{amssymb}
\usepackage{amsfonts}
\usepackage{epsfig}
\usepackage{amsmath}
\usepackage{amscd}
\usepackage{psfrag}
\usepackage[all]{xy}

\newtheorem{defn0}{Definition}[section]
\newtheorem{prop0}[defn0]{Proposition}
\newtheorem{thm0}[defn0]{Theorem}
\newtheorem{lemma0}[defn0]{Lemma}
\newtheorem{corollary0}[defn0]{Corollary}
\newtheorem{example0}[defn0]{Example}
\newtheorem{conjecture0}[defn0]{Conjecture}
\newtheorem{notation0}[defn0]{Notation}

\theoremstyle{remark}
\newtheorem{remark0}[defn0]{Remark}

\newenvironment{prop}{\begin{prop0}}{\end{prop0}}
\newenvironment{thm}{\begin{thm0}}{\end{thm0}}
\newenvironment{lem}{\begin{lemma0}}{\end{lemma0}}
\newenvironment{cor}{\begin{corollary0}}{\end{corollary0}}
\newenvironment{example}{\begin{example0}\rm}{\end{example0}}
\newenvironment{rem}{\begin{remark0}\rm}{\end{remark0}}

\renewenvironment{proof}{\noindent {\textsc{Proof.}}}{$\square$ \vspace{3mm}}

\newcommand{\K}{{\mathcal K}}
\newcommand{\Q}{{\mathcal Q}}
\newcommand{\E}{{\mathcal E}}

\newcommand{\T}{{\mathcal T}}

\newcommand{\LL}{{\mathcal L}}
\newcommand{\II} { \textbf{I}_R }

\newif\ifprivate
\privatetrue

\def\???{\ifprivate {\bf {???}} \marginpar{{\Huge {\bf ?}}}
\else \fi}

 \DeclareMathOperator{\mult}{mult}
 \DeclareMathOperator{\lcm}{lcm}

 \DeclareMathOperator{\sC}{\scriptscriptstyle \textit{C}}

 \numberwithin{equation}{section}

\begin{document}

\title{On curves on sandwiched surface singularities}

\author{Jes\'{u}s Fern\'{a}ndez-S\'{a}nchez}
\thanks{Partially supported by
CAICYT BFM2002-012040, Generalitat de Catalunya 2001SGR 00071 and EAGER,
European Union Contract HPRN-CT-2000-00099.
\newline{AMS 2000 subject classification 14H20; 14H50; 13B22; 14B05; 14E05}}

\address{Departament de Matem\`atica Aplicada I \\ Universitat Polit\`ecnica de Catalunya \\
Av. Diagonal 647, 08028-Barcelona, Spain.}
 \email{jesus.fernandez.sanchez@upc.edu}

\vspace{6mm} \begin{abstract} Fixed a point $O$ on a
non-singular  surface $S$ and a complete
$\mathfrak{m}_O$-primary ideal $I$ in ${\mathcal O}_{S,O}$,
the curves on the surface $X=Bl_I(S)$ obtained by
blowing-up $I$ are studied in terms of the base points of
$I$. Criteria for the principality of these curves are
obtained. New formulas for their multiplicity, intersection
numbers and order of singularity at the singularities of
$X$ are given. The semigroup of branches going through a
sandwiched singularity is effectively determined, too.
\end{abstract}

\maketitle

\section{Introduction}
Sandwiched singularities are normal surface singularities
which birationally dominate a non-singular surface. They
are rational surface singularities and among them are
included all cyclic quotients and minimal surface
singularities. The original interest in sandwiched
singularities comes from a question posed by Nash in the
early sixties: \emph{does a finite succession of Nash
transformations or normalized Nash transformations resolve
the singularities of a reduced algebraic variety?} Hironaka
had proved in \cite{hironaka} that after a finite
succession of normalized Nash transformations one obtains a
surface $X$ having only sandwiched singularities. Some
years later, in \cite{Spivak}, M. Spivakovsky proves that
sandwiched singularities are resolved by normalized Nash
transformations, thus giving a positive answer to the
original question posed by Nash for the case of surfaces.
Since then, a constant interest in sandwiched singularities
has been shown, and they have been deeply studied as a nice
testing ground for the Nash and the wedge problem by
Lejeune-Jalabert and Reguera in \cite{LR}. A recent paper
by Reguera proves the Nash conjecture for sandwiched
singularities \cite{Reguera04}. Sandwiched singularities
have been also studied from the point  of view of
deformation theory by de Jong, van Straten in \cite{JS} and
by Gustavsen in \cite{stolen}.

\vspace{3mm}
Any sandwiched surface singularity can be obtained by
blowing up a complete $\mathfrak{m}_O$-primary ideal $I$ in
a local ring ${\mathcal O}_{S,O}$ of a non-singular surface
$S$. The main purpose of this  paper is to study the curves
(effective Weil divisors) on the sandwiched surface
$X=Bl_I(S)$ through their relationship to the infinitely
near base points of $I$.
One motivation for this study comes from \cite{Lip69} and
the study of the divisor class group for rational surface
singularities in general. Another motivation comes from the
fact that the study of the curves on a surface singularity
may lead to a deeper understanding of the singularity
itself (see \cite{Nash}). Criteria for the principality of
Weil divisors on $X$ are given in this paper. A number of
invariants for the curves on $X$ (such as their
multiplicity and their order of singularity) are also
recovered in terms of the relative position of their
projection to $S$ with the base points of $I$.
We keep the point of view of
\cite{moi1} and make use of the theory of infinitely near
points as revised and developed by Casas-Alvero in
\cite{CasasBook}. Relevant to our purposes is the fact that
any complete ideal has a cluster of infinitely near base
points that in turn determines the ideal.

The organization of the paper is as follows. Concepts and
well-known facts about infinitely near points and
sandwiched surface singularities are reviewed from
\cite{CasasBook} and \cite{moi1} in section 2. We also give
a formula for the multiplicity of curves on $X=Bl_I(S)$
going through a sandwiched singularity in terms of the base
points of the ideal $I$ (Corollary \ref{multip_curve}).
Section 3 deals with the existence of equations for curves
on a sandwiched surface. The main result of the paper is
Theorem \ref{Cartier}, which provides different criteria
for their principality. In Section 4 we derive easy
consequences and show some examples to illustrate these
results. In Section 5, we suggest a procedure to compute
effective Cartier divisors containing a given (not
necessarily principal) curve $C$. To this aim, we introduce
a slight modification of the unloading procedure called
\emph{partial unloading}. It is explained in the Appendix
and it may be useful in any other context where precise
control of unloading is needed. The procedure of Section 5
gives rise to flags of clusters depending on the ideal $I$
(and hence, on the surface $X$) and the curve $C$. These
flags keep deep information about the strict transform of
$C$ and they are used in Section 6 to infer formulas for
its order of singularity (Proposition \ref{delta/general})
and to compute effectively the semigroup of a branch on any
singularity of $X$ (Proposition \ref{semlemma}).

\vspace{2mm}
\section{Preliminaries}

Throughout this work the base field is the field
$\mathbb{C}$ of complex numbers. A standard reference for
most of the material treated here is the book by
Casas-Alvero \cite{CasasBook}. Let $(R,{\mathfrak{m}_O})$
be a regular local two-dimensional $\mathbb{C}$-algebra and
$S=Spec(R)$. A \emph{cluster} of points of $S$ with origin
$O$ is a finite set $K$ of points infinitely near or equal
to $O$ such that for any $p\in K$, $K$ contains all points
preceding $p$. By assigning integral multiplicities
$\nu=\{\nu_p\}$ to the points of $K$, we get a
\emph{weighted cluster} ${\mathcal K}=(K,\nu)$, the
multiplicities $\nu$ being called the \emph{virtual
multiplicities} of $\K$. We write $p\rightarrow q$ if $p$
is proximate to $q$. If
\[\rho_p=\nu_p-
\sum_{q \rightarrow p} {\nu _q }
\]
is the \emph{excess} at $p$ of ${\mathcal K}$,
\emph{consistent} clusters are those clusters with no
negative excesses. Moreover, $\K_+=\{p\in K|\rho_p>0\}$ is
the set of \emph{dicritical} points of $\K$. \emph{Strictly
consistent} clusters are consistent clusters with no points
of virtual multiplicity zero.

If ${\mathcal K}$ is a weighted cluster, the equations of
all curves going through it define a complete
$\mathfrak{m}_O$-primary ideal $H_{\mathcal K}$ in $R$ (see
\cite{CasasBook} 8.3). Two weighted clusters ${\mathcal K}$
and ${\mathcal K}'$ are \emph{equivalent} if $H_{\mathcal
K}=H_{{\mathcal K}'}$. Any complete
$\mathfrak{m}_O-$primary ideal $J$ in $R$ has a cluster of
base points, denoted by  $BP(J)$, which consists of the
points shared by, and the multiplicities of, the curves
defined by a generic element of $J$. Moreover, the maps
$J\mapsto BP(J)$ and ${\mathcal K}\mapsto H_{\mathcal K}$
are reciprocal isomorphisms between the semigroup $\II$ of
complete $\mathfrak{m}_O$-primary ideals in $R$ and the
semigroup of strictly consistent clusters (see
\cite{CasasBook} 8.4.11 for details). If $p\in {\mathcal
N}_O$, we denote by $I_p$ the simple ideal generated by the
equations of the branches going through $p$ and by
${\mathcal K}(p)$ the weighted cluster corresponding to it
by the above isomorphism. Clearly, $I_p=H_{\K(p)}$. In this
language, if $I\in \II$ and $\K=BP(I)$ is the cluster of
base points of $I$, we have that $\K=\sum_{i=1}^n\alpha_i
\K(p_i)$ if and only if
$I=\prod_{i=1}^nI_{p_i}^{\alpha_{i}}$ is the Zariski
factorization of $I$ into simple ideals. Hence, the
exponent $\alpha_p$ of $I_p$ in the factorization of $I$
equals the excess of $\K$ at $p$, so we can write
\begin{eqnarray}\label{for:factorization}
I=\prod_{p\in \K_+}I_p^{\rho_p}.
\end{eqnarray}
If $\K=(K,\nu)$ is consistent, then the \emph{virtual
codimension} of $\K$ equals (Proposition 4.7.1 of
\cite{CasasBook})
\begin{eqnarray}
\label{for:cod_cluster}c(\K)=\sum_{p\in
K}\frac{\nu_p(\nu_p+1)}{2}}={\dim_{\mathbb{C}}(\frac{R}{H_{\K}}).
\end{eqnarray}

Consistent clusters are characterized
as those clusters whose virtual multiplicities may be realized effectively by
some curve on $S$. If ${\mathcal K}$ is not consistent, $\widetilde{\mathcal
K}$ is the cluster given rise from ${\mathcal K}$ by {\em unloading}.
Equivalently, $\widetilde{\mathcal K}$ is the unique consistent cluster which
is equivalent to ${\mathcal K}$ and has the same points (see \cite{CasasBook}
\S 4.2 and \S 4.6 for details).

\vspace{3mm} If $\pi_K:S_K\longrightarrow S$ is the
composition of the blowing-ups of all points in $K$, write
$\E_K$ for the exceptional divisor of $\pi_K$ and
$\{E_p\}_{p\in K}$ for its irreducible components. Use
$|\,\,{\cdot}\,\,|$ as meaning intersection number and
$[\,\,,\,\,]_P$ as intersection multiplicity at $P$. We
denote by $\mathbf{A}_K=-\textbf{P}_K^t\textbf{P}_K$ the
\emph{intersection matrix} of $\E_K$, where $\textbf{P}_K$
is the \emph{proximity matrix} of $K$. If $C$ is a curve on
$S$, $v_p(C)$ is the value of $C$ relative to the
divisorial valuation associated with $E_p$. Then,
$v_p(C)=e_p(C)+\sum_{p\mapsto q}v_q(C)$. We write
$\mathbf{v}_{\K}=(v_p^{\K})_{p\in K}$ for the vector of
virtual values of $\K$. Moreover, if $\widetilde{C}^{K}$ is
the strict transform of $C$ on $S_K$, we have
\begin{eqnarray}\label{for:int/comp}
|\widetilde{C}^{K}\cdot E_p|_{S_K}=e_p(C)-\sum_{q\in K,q \mapsto p}e_q(C), \mbox{ for all
$p\in K$.}
\end{eqnarray}
and we have the equality (\emph{projection formula}) for
$\pi_K$: $|\pi_K^*(C)\cdot B|_{S_K}=[C,(\pi_K)_*(B)]_O$ for
any curve $B$ on $S_K$. We will make use also of the
Noether formula for the intersection multiplicity of two
curves at $O$ (Theorem 3.3.1 of \cite{CasasBook}):
\begin{eqnarray}\label{for:Noether}
[C,C']_O=\sum_{p}e_p(C)e_p(C').
\end{eqnarray}
For consistent clusters, write $[{{\mathcal
K},C]_O=\sum_{p\in K}\nu_p e_p(C)}$ for the
\emph{intersection multiplicity} of ${\mathcal K}$ and a
curve $C$ on $S$, and ${\mathcal K}^2=\sum_{p\in K}\nu_p^2$
for the \emph{self-intersection} of ${\mathcal K}$.

\vspace{3mm} We will also make use of the following results
for rational surface singularities. Let $(X,Q)$ be a
rational surface singularity and let $g:X'\rightarrow X$ be
any resolution of it. Write $\{E_i\}_{i=1,\ldots,n}$ for
the irreducible components of the exceptional locus of $g$.
If $A$ is a curve on $X$ (i.e. an effective Weil divisor),
it is possible to associate to $A$ a $\mathbb{Q}$-Cartier
exceptional divisor $D_A$ on $X'$, being defined by the
condition (see II (b) of \cite{Mumf61})
\begin{eqnarray*} \label{eqmunf}
|D_A\cdot E|_{X'} \, = \, -|\widetilde{A}^{X'}\cdot E|_{X'}
\quad \forall i\in \{1,\ldots,n\}.
\end{eqnarray*}
On the other hand,  if $D$ is a divisor on $X'$ such that
$|D\cdot E_i|=0$ for each exceptional component, then there
exists an element $h\in \mathfrak{m}_{X,Q}$ such that $D$
is the total transform of the curve defined by $h=0$ (see
the proof of Theorem 4 of \cite{Art66}). Because of this,
we have that
\begin{lem}
\label{mumf} $A$ is Cartier if and only if $D_A$ is a divisor on $X'$.
\end{lem}

\vspace{2mm} \subsection{Sandwiched surface singularities}

In this section, we fix notation and recall some facts
concerning sandwiched singularities and base points of
ideals which will be useful throughout this paper. The main
references are \cite{Spivak} and \cite{moi1}. A
\emph{curve} will be an effective Weil divisor on a
surface.

Let $I\in \II$. The ideal $I$ is fixed throughout the paper, so that no confusion should
arise if the notation introduced from now on does not reflect its dependence of $I$. Write
$\K=(K,\nu)$ for the cluster of base points of $I$ and $\pi_I:~X=Bl_I(S)\longrightarrow
S$ for the blowing-up of $I$. The singularities of $X$ are by definition sandwiched
singularities and we have a commutative diagram
\[\xymatrix{  {S_K} \ar[r]^{f}\ar[rd]_{\pi_K} & {X} \ar[d]^{\pi_I} \\ &  {S}}\]
where the morphism $f$, given by the universal property of the blowing-up, is the minimal
resolution of the singularities of $X$ (Remark 1.4 of \cite{Spivak}).

If $Q$ is in the exceptional locus of $X$, write
$\mathfrak{m}_Q$ for the maximal ideal in ${\mathcal
O}_{X,Q}$ and ${\mathcal M}_Q$ for the ideal sheaf on $X$
associated with $Q$. Then, the ideal $I_Q:=\pi_*({\mathcal
M}_QI{\mathcal O}_X)\subset I$ is complete,
$\mathfrak{m}_O$-primary and has codimension one in $I$ and
moreover, the map $Q\mapsto I_Q$ defines a bijection
between the set of points in the exceptional locus of $X$
and the set of complete $\mathfrak{m}_O$-primary of
codimension one in $I$ (Theorem 3.5 of \cite{moi1}).

It can be seen that if $q \in E_K$, and ${\mathcal K}_q$ is
the weighted cluster obtained from $\K$ by adding $q$ as a
simple point, then $H_{\K_q}\subset I$ has codimension one
in $I$, and every complete $\mathfrak{m}_O$-primary ideal
of codimension one in $I$ has this form (Lemma 3.1 of
\cite{Edinburgh}). Therefore, for every point $Q$ in the
exceptional locus of $X$, there exists some $q\in E_K$ such
that $I_Q=H_{\K_q}$
It follows also that any $q\in E_K$ can be mapped to the
point of $X$ corresponding to $H_{{\mathcal K}_q}\subset I$
by the above bijection. This correspondence is surjective,
but not injective since the pre-image by this map of any
singular point of $X$ consist of the points lying on the
exceptional components contracting by $f$ to that
singularity. On the other hand, non-singular points of $X$
are characterized as those points for which there is only
one $q\in E_K$ such that $I_Q=H_{\K_q}$ (Proposition 4.4 of
\cite{moi1}). In any case, we denote by $\K_Q$ the cluster
of base points of $I_Q$, $\K_Q=BP(I_Q)$.

For any sandwiched singularity $Q\in X$, write $T_Q=\{p\in
K\mid f_*(E_p)=Q\}$ so that $\{E_p\}_{p\in T_Q}$ is the set
of the exceptional components of $S_K$ contracting to $Q$.
In particular, we have that
\begin{eqnarray}
\label{for:excsing} \{p\in K \mid \rho_p=0\}=\bigcup_{Q\in
Sing(X)}T_Q.
\end{eqnarray}

Recall that there is a bijection between the set of simple
ideals $\{I_p\}_{p\in \K_+}$ in the (Zariski) factorization
of $I$ (see (\ref{for:factorization})) and the set of
irreducible components of $\pi_I^{-1}(O)$ (\cite{Spivak}
Corollary 1.5; see also \cite{Lip69} Proposition 21.3). We
write $\{L_p\}_{p\in \K_+}$ for the set of these
components. If $J\in \II$, we write $\LL_J=\sum_{p\in
\K_+}v_p(J)L_p$ and if $C$ is a curve on $S$,
$\LL_C=\sum_{p\in \K_+}v_p(C)L_p$. For $p\in \K_+$, write
${\mathcal L}_p=\sum_{q\in \K_+}v_q(I_p)L_q$.

The following proposition is easy, but will be very useful
in the future:

\begin{prop} [Projection formula for $\pi_I$]
If $C,D$ are curves on $S$, we have that
\[[C,D]_O=|C^*\cdot \widetilde{D}|_X.\]Also, $|C^*\cdot
L_p|_X=0$, for any $p\in \K_+$.
\end{prop}

\begin{proof}
Let $g:X'\rightarrow X$ be a resolution of $X$. To prove
the claim it is enough to apply the projection formula for
$\pi_K$ and for $g$.
\end{proof}

The following technical result is a generalization of Lemma
4.2 of \cite{moi2}. Given two points $q,p\in K$, the chain
$ch(q,p)$ is the subgraph of the dual graph $\Gamma_K$ of
$E_K$ of all vertices and edges between the vertices
representing $q$ and $p$ in $\Gamma_K$ (see \S 4.4 of
\cite{CasasBook}).

\vspace{3mm}
\begin{lem} \label{new_multiplicities} Let $Q\in X$ be singular and write $\K_Q=(K,\nu^{(Q)})$.
\begin{itemize}
\item[(a)] By taking the partial order relation of \emph{being infinitely near
to}, there exists a unique minimal point in $T_Q$. We will denote this point by
$O_Q$;
 \item[(b)] $\nu^{(Q)}_{O_Q}=\nu_{O_Q}+1$. If $p\in K$ and $p\neq O_Q$, then
 $\nu_p-1 \leq \nu^{(Q)}_p \leq \nu_p$.

\end{itemize}
\end{lem}

\begin{proof}
(a) Assume that $q_1,q_2\in T_Q$ are different and minimal
among the points of $T_Q$. Let $u_0\in K$ be such that
$q_1$ and $q_2$ are infinitely near to $u_0$ and maximal
with this property. Then, $u_0\in ch(q_1,q_2)$. By the
minimality of $q_1$ and $q_2$, $u_0\notin T_Q$. By
Zariski's Main Theorem (see, for example, Theorem V 5.2 of
\cite{Hartshorne}), the union of the components
$\{E_p\}_{p\in T_Q}$ is connected and by Theorem 1.7 of
\cite{Brieskorn}, $\E_K$ contains no cycles. This leads to
contradiction, and (a) follows.

To prove (b), denote by $S^{(Q)}$ the surface obtained by
blowing up the points preceding $O_Q$, so that $O_Q$ is a
proper point of $S^{(Q)}$ ($S^{(Q)}=S$ if $O_Q=O$). Denote
by $\varphi_{O_Q}:R\longrightarrow {\mathcal
O}_{S^{(Q)},O_Q}$ the morphism induced by blowing up. By
Corollary 3.6 of \cite{moi1}, we know that
$v_p(I_Q)=v_p(I)$ and so, $\nu^{(Q)}_p=\nu_p$ for each $p$
preceding $O_Q$. Therefore, the exceptional component of
the total transform on $S^{(Q)}$ of curves going sharply
through $\K$ or $\K_Q$ are equal. Write $z$ for an equation
of this exceptional component at $O_Q$. Then, the ideal
$\breve{I}_Q$ (resp. $\breve{I}$) generated by
$z^{-1}\varphi_{O_Q}(I_Q)$ (resp. $z^{-1}\varphi_{O_Q}(I)$)
is complete and $\mathfrak{m}_{O_Q}$-primary in ${\mathcal
O}_{S^{(Q)},O_Q}$. Since $I_Q\subset I$, we have
$\breve{I}_Q\subset \breve{I}$. Moreover, the base points
of $\breve{I}$ (resp. $\breve{I}_Q$) equal the base points
of $I$ (resp. $I_Q$) infinitely near or equal to $O_Q$.
Direct computation using (\ref{for:cod_cluster}) shows that
$\dim_{\mathbb{C}}(\frac{\breve{I}}{\breve{I}_Q})=1$. Now,
the proof of (b) follows as in Lemma 4.2 of \cite{moi2}.
The details are left to the reader.
\end{proof}

\noindent The above lemma shows that the virtual
multiplicities of $\K_Q$ differ as much in one from the
original multiplicities of $\K$. The set \[{\mathcal
B}_Q=\{p\in K \mid \nu^{(Q)}_p=\nu_p-1\}\] will play an
important role in the sequel. Given $p\in K$, we denote
$\varepsilon_p=-1$ if $p\in {\mathcal B}_Q$,
$\varepsilon_{O_Q}=1$ and $\varepsilon_p=0$, otherwise.

\vspace{2mm} Let us quote some consequences of
\ref{new_multiplicities}. The first one follows immediately
from $\dim_\mathbb{{C}}\frac{I}{I_Q}=1$ by using
(\ref{for:cod_cluster}).

\begin{cor}
\label{equilBQ} $\displaystyle \nu_{O_Q}=\sum_{p\in
{\mathcal B}_Q}\nu_p$.
\end{cor}

We also obtain a new formula for the multiplicity at $Q$ of
curves on $X$ in terms of $O_Q$ and the points of
${\mathcal B}_Q$. Namely,

\begin{cor} \label{multip_curve}
Let $Q$ be a point in the exceptional locus of $X$ and let $C$ be a curve on $S$. Then
$$\mult_Q(\widetilde{C})=e_{O_Q}(C)-\sum_{p\in {\mathcal B}_Q}e_p(C).$$
\end{cor}

\begin{proof}
First of all, notice that if $Q$ is not singular,  then
${\mathcal B}_Q$ is empty and $O_Q=O$, so there is nothing
to prove. Hence, we assume that $Q$ is singular. By the
projection formula applied to $f$, we have that
\begin{eqnarray}
\label{for:mult_curve} \mult_Q(\widetilde{C}) & = &
|\widetilde{C}^{K}\cdot Z_Q|_{S_K}=|\widetilde{C}^{K}\cdot
(\E_{I_Q}-\E_I)|_{S_K}
\end{eqnarray}
the last equality by Corollary 3.6 of \cite{moi1}. Now,
since the sheaves $I_Q {\mathcal O}_{S_K}$ and $I{\mathcal
O}_{S_K}$ are invertible, we can take two curves $B$ and
$B'$ going sharply through $\K$ and $\K_Q$, respectively
and such that $\widetilde{C}_I$ and $\widetilde{C}_{I_Q}$
share no points with $\widetilde{C}^{K}$ on $S_K$. Then, by
the projection formula applied to $\pi_K$, $[C,B]_O =
|\widetilde{C}^{K}\cdot (\widetilde{C}_I^{K}+\E_I)|_{S_K}=
|\widetilde{C}^{K}\cdot \E_I|_{S_K}$ and similarly, $
[C,B']_O = |\widetilde{C}^{K}\cdot \E_{I_Q}|_{S_K}$. Hence,
by (\ref{for:mult_curve}) above, $ \mult_Q(\widetilde{C}) =
[C,B']_O -[C,B]_O$, and by Noether's formula (see
(\ref{for:Noether})), we have
$\mult_Q(\widetilde{C})=\sum_p e_p(C)(\nu^{(Q)}_p-\nu_p)$.
Finally, by (b) of \ref{new_multiplicities}, we obtain
$\mult_Q(\widetilde{C})=e_{O_Q}(C)-\sum_{p\in {\mathcal
B}_Q}e_p(C)$ as claimed.
\end{proof}

An easy formula for the multiplicity of a sandwiched
singularity $Q$ is also obtained in terms of the set
${\mathcal B}_Q$. Moreover, it shows that the number of
these points is an invariant of the singularity and so, it
is independent of the particular ideal blown-up to obtain
it. Recall that $R={\mathcal O}_{S,O}$ is said to be
\emph{maximally regular} in ${\mathcal O}_{X,Q}$ if there
are no regular rings $R_0\neq R$ such that $R\subset R_0
\varsubsetneq {\mathcal O}_{X,Q}$ (see \cite{HunSally}).

\begin{cor}\label{multsand} Let $Q$ be a point in the exceptional locus of $X$.
Then,
\begin{itemize}
\item[(a)] $\mult_Q(X)=1+\sharp {\mathcal B}_Q$;
 \item[(b)] $\mult_Q(X)\leq 1+\nu_{O}$, and the equality holds if and only if $R$ is maximally regular in ${\mathcal
O}_{X,Q}$ and all the points in ${\mathcal B}_Q$ are simple.
\end{itemize}
\end{cor}

\begin{proof}
It is enough to compute the multiplicity of a transverse
hypersurface section of $(X,Q)$, that is, of the strict
transform of a curve $C$ going sharply through $\K_Q$, see
Theorem 3.5 of \cite{moi1}. (a) follows from
\ref{multip_curve} and \ref{equilBQ} by direct computation.
Now, the inequality in (b) is obvious if $Q$ is regular, so
we assume that $Q$ is singular. Clearly, $\sharp {\mathcal
B}_Q\leq \sum_{p\in {\mathcal B}_Q}\nu_p$. Again by
\ref{equilBQ}, we have
\begin{eqnarray*}
\mult_Q(X) \leq 1+\sum_{p\in {\mathcal B}_Q}\nu_p= 1+\nu_{O_Q} \leq 1+\nu_O,
\end{eqnarray*}
and the equality holds if and only if the virtual
multiplicity of all the points in ${\mathcal B}_Q$ is one
and $\nu_O=\nu_{O_Q}$. Now, if $O_Q\neq O$, there must be
some $q$ preceding $O_Q$ with $\rho^{\K_Q}_q>0$ (recall
that $O_Q$ is minimal in $T_Q$, see
\ref{new_multiplicities}) and so, $\nu_q>\nu_{O_Q}$. Hence,
$\nu_O \geq \nu_q > \nu_{O_Q}$ and $\mult_Q(X)<1+\nu_O$.
Now, this is equivalent to  the maximally regular
condition, since ${\mathcal O}_{X,Q}$ can always be
projected birationally into the surface obtained by blowing
up all the points preceding $O_Q$, i.e. the surface where
$O_Q$ is lying as a proper point, so $R\subset {\mathcal
O}_{O_Q}\subsetneq {\mathcal O}_{X,Q}$. This completes the
proof.
\end{proof}

\section{Criteria for the principality of curves through a sandwiched singularity}

\begin{thm}
\label{Cartier} Let $C$ be a curve on $S$ and write
$\LL_C=\sum_{u\in \K_+}v_u(C)L_u$. The following conditions
are equivalent:
\begin{itemize}
\item[(i)] $\widetilde{C}$ is a Cartier divisor on $X$;
 \item[(ii)] $\LL_C \in \bigoplus_{u\in \K_+} \mathbb{Z} \mathcal{L}_u$;
 \item[(iii)] there exists a
 curve $B$ on $S$ such that $\LL_{B}=\LL_C$ and $\widetilde{B}$ goes through no
singularity of $X$;
 \item[(iv)] if $J_C=\{g\in R \mid v_p(g)\geq v_p(C), \forall p \in \K_+\}$
 and ${\mathcal Q}_{\sC}=BP(J_C)$, then every dicritical point of ${\mathcal Q}_{\sC}$ is a
dicritical point of $\K$.
\end{itemize}
\end{thm}

\vspace{2mm} Before proving \ref{Cartier} we need a
technical result.


\begin{lem}
\label{basis}
\begin{itemize}
\item[(a)] The set $\{\E_{I_p}\}_{p\in K}$ is a basis of
the $\mathbb{Q}$-vector space generated by $\{E_p\}_{p\in
K}$. The matrix of the change of basis from
$\{\E_{I_p}\}_{p\in K}$ to $\{E_p\}_{p\in K}$ is
$-\mathbf{A}_K$.

\item[(b)] The set $\{\mathcal{L}_u\}_{u\in \K_+}$ is
a basis of the $\mathbb{Q}$-vector space generated by
$\{L_u\}_{u\in \K_+}$.
\end{itemize}
\end{lem}

\begin{proof}
We have $\mathbf{A}_K(v_q(I_p))_{q\in K}=-\mathbf{1}_p$
where $\mathbf{1}_p$ is the $K$-vector having all its
entries equal to $0$ but the corresponding to $p$ which is
$1$. Thus, $\mathbf{A}_K \E_{I_p}=-E_p$ and the claim (a)
follows. To prove (b), it is enough to show that the
$\{\mathcal{L}_u\}_{u\in \K_+}$ are linearly independent.
Assume that there exist rational numbers $\{a_u\}_{u\in
\K_+}$ such that
\begin{eqnarray}
\label{for:auxC} \sum_{u\in \K_+}a_u \mathcal{L}_u=0.
\end{eqnarray}
By multiplying by an integer, we may assume that the
$a_u\in \mathbb{Z}$ for all $u\in \K_+$. Now, for each
$u\in \K_+$ take $\gamma_u$ a curve going sharply through
${\mathcal K}(u)$ and missing all points after $u$ in $K$,
and write
$$\sum_{u\in \K_+}a_u \gamma_u=C_1-C_2$$ where $C_1=\sum_{a_u>0}a_u
\gamma_u$ and $C_2=\sum_{a_u<0}(-a_u) \gamma_u$. Then, by
(\ref{for:auxC}) we have that $\LL_{C_1}=\sum_{a_u>0}a_u
\mathcal{L}_u $ and $\LL_{C_2}=\sum_{a_u<0}(-a_u)
\mathcal{L}_u $ are equal. Hence, by taking total
transforms on $S_K$, we see that
\begin{eqnarray*}
\sum_{a_u>0}a_u \E_{I_u}=\sum_{a_u<0}(-a_u) \E_{I_u}
\end{eqnarray*} against (a).
\end{proof}

\vspace{2mm} \noindent \textsc{Proof of \ref{Cartier}}. We prove that $(i)\Rightarrow
(ii)\Rightarrow (iii)\Rightarrow (iv)\Rightarrow (i)$.

\vspace{2mm} \noindent$(i)\Rightarrow (ii)$. Assume that $\widetilde{C}$ is Cartier on
$X$ and hence, that $\LL_C$ is so (as $C^*=\widetilde{C}+\LL_C$ is
always Cartier). Then, by \ref{mumf}, $f^*(C)$ is a divisor on $S_K$ and so, the
coefficients of the components $\{E_p\}_{p\in T_Q}$ in $f^*(\LL_C)$ are integers, i.e.
\begin{eqnarray}
\label{for:coefL} f^*(\LL_C)=\sum_{q\in K}b_q E_q \quad
\mbox{with $b_q\in \mathbb{Z}$.}
\end{eqnarray}
On the other hand, by (b) of \ref{basis}, we can write
$\LL_C=\sum_{u\in \K_+}a_u \mathcal{L}_u$, with $\{a_u\}$
rational numbers. Now, if $u\in \K_+$, then
$f^*(\mathcal{L}_u)=\E_{I_u}$ and $f^*(\LL_C)=\sum_{u\in
\K_+}a_u \E_{I_u}$ is the expression of $f^*(\LL_C)$ in the
basis $\{\E_{I_p}\}_{p\in K}$. Therefore, by (a) of
\ref{basis} and the equality (\ref{for:coefL}),
$(a_u)_{u\in K}=-\mathbf{A}_K (b_u)_{u\in K}$ and hence all
the $a_u$ are integers.

\vspace{2mm} \noindent $(ii)\Rightarrow (iii)$. Assume that
$\LL_C = \sum_{u\in \K_+}a_u \mathcal{L}_u,$ with $a_u\in
\mathbb{Z}$. If $u\in \K_+$, the projection formula applied
to $\pi$ gives that $|C^*\cdot
L_u|_X=|(\widetilde{C}+\LL_C)\cdot L_u|_X=0$. Hence,
\begin{eqnarray}
|\LL_C\cdot L_u|_X= -|\widetilde{C}\cdot L_u|_X.
\end{eqnarray}
In particular, if $\gamma_p$ is a curve going sharply
through $\K(p)$ ($p\in \K_+$) and missing all points in $K$
after $p$, then
\begin{eqnarray*}
|\mathcal{L}_p\cdot L_u|_X=-|\widetilde{\gamma_p}\cdot L_u|_X=\left \{
\begin{array} {ll} -1 &  \mbox{ if }p=u
\\ 0 &  \mbox{ otherwise.} \end{array} \right.
\end{eqnarray*} It follows that
\[|\widetilde{C}\cdot L_u|_X=-|\LL_C\cdot L_u|_X=-\sum_{p\in
\K_+}a_p|\mathcal{L}_p\cdot L_u|_X=a_u\] and thus, $a_u\geq
0$ for all $u\in \K_+$. Now, let $B$ be a curve going
sharply through $\T=\sum_{u\in \K_+}{a_u}{\K(u)}$ and
missing those points of $K$ not contained in the underlying
cluster of $\T$ (Corollary 4.2.8 of \cite{CasasBook}). The
strict transform $\widetilde{B}\subset X$ cuts
transversally each exceptional component $L_u$ at $a_u$
different points and goes through no singularity of $X$.
Moreover, $\LL_{B}=\sum_{u\in \K_+} a_u
\mathcal{L}_u=\LL_C$.

\vspace{2mm} \noindent $(iii)\Rightarrow (iv)$. Since
$\LL_{B}=\LL_C$, we have $v_p(B)=v_p(C)=v_p(J_C)$ for all
$p\in \K_+$, and hence, $B$ is defined by an element of
$J_C$. Therefore, $v_p(B)\geq v_p(J_C)$ for all $p$, and so
\begin{eqnarray}
\label{auxxCartier} \E_{B}\geq \E_{J_C}=\sum_{p\in
K}v_p(J_C)E_p.
\end{eqnarray}
On the other hand, since $\widetilde{B}$ goes through no
singularities of $X$, the total transform of
$\widetilde{B}$ by $f$ has no exceptional part, and
$\pi_K^*(B)=f^*(\widetilde{B}+\LL_C)=\widetilde{B}^{K}+f^*(\LL_C)$.
Thus, $\E_{B}=f^*(\LL_C)$. Now, if $C_1$ goes sharply
through ${\mathcal Q}_{\sC}$, $\E_{C_1}=\E_{{\mathcal
Q}_{\sC}}$ and $\LL_{C_1}=\LL_C$. Assume that there exists
$q\in K\setminus \K_+$ so that ${\mathcal Q}_{\sC}$ has
positive excess at it. Then $\rho_q=0$ and by the equality
(\ref{for:excsing}) in page \pageref{for:excsing}, there is
some singularity $Q$ in $X$ such that $q\in T_Q$. Moreover,
$\widetilde{C_1}$ goes through $Q$ and
$f^*(\widetilde{C_1})=\widetilde{C_1}^{K}+D_{\widetilde{C_1}}$
where $D_{\widetilde{C_1}}=\sum_{u\in K\setminus \K_+}b_u
E_u$ is a \emph{non-zero effective} divisor. Thus,
\begin{eqnarray*}
\E_{C_1} & = & D_{C_1}+f^*(\LL_{C_1})=
D_{C_1}+f^*(\LL_C)=D_{C_1}+\E_{B}>\E_{B}
\end{eqnarray*} against (\ref{auxxCartier}). Therefore,
${\mathcal Q}_{\sC}$ has excess $0$ at every $q\in
K\setminus \K_+$.

\vspace{2mm} \noindent $(iv)\Rightarrow (i)$. Assume that
every dicritical point of ${\mathcal Q}_{\sC}$ is also
dicritical for $\K$. Let  $B$ be a curve going sharply
through ${\mathcal Q}_{\sC}$ and such that no points in
$K\setminus {\mathcal Q}_{\sC}$ belong to it. Then,
$\widetilde{B}^{K}$ intersects no components $E_p$ with
$p\in K\setminus \K_+$. Since the direct image of
$\widetilde{B}^{K}$ by $f$ is the strict transform
$\widetilde{B}$  on $X$, it follows that $\widetilde{B}$
goes through no singularities of $X$ and is Cartier. Using
that $\LL_{B}=\LL_C$ and that $B^*=\widetilde{B}+\LL_{B}$
and $C^*=\widetilde{C}+\LL_C$ are Cartier divisors, we
deduce that $\LL_C$ and $\widetilde{C}$ are so. \qed

\vspace{3mm} \emph{Primitive} singularities are those
singularities that can be obtained by blowing-up a simple
complete ideal (Definition I.3.1 of \cite{Spivak}). In this
case, \ref{Cartier}  has a very easy formulation.

\begin{cor}
Let $I_p\subset R$ be a simple ideal and $X=Bl_{I_p}(S)$. If $C\subset S$ is a
curve, $\widetilde{C}$ is a Cartier divisor on $X$ if and only if $v_p(C)$ is a
multiple of ${\mathcal K}(p)^2$. Moreover,
\[m_{\widetilde{C}}=\frac{\lcm (v_p(C),\K(p)^2)}{v_p(C)}\] is the minimal integer
$m$ such that $m \widetilde{C}$ is a Cartier divisor.
\end{cor}

\begin{proof}
It is clear that $\LL_C=v_p(C)L_p$ and $\mathcal{L}_p=v_p(I_p)L_p$. Since
$\K(p)^2=v_p(I_p)$, we infer that $\LL_C\in \mathbb{Z}\mathcal{L}_p$ if and only if
$v_p(C)\in ({\mathcal K}(p)^2)$. The second claim follows by direct computation.
\end{proof}

Clearly, the question of whether $\widetilde{C}$ on $X$ is
Cartier is local as it depends only on the existence of an
equation for $\widetilde{C}$ near the singularities of $X$
(recall that every Weil divisor on a neighbourhood of a
regular point is principal, see \cite{Hartshorne} II.6.11).
If $Q$ is in the exceptional locus of $X$, we write
$\K_+^Q$ for the set of points $p\in \K_+$ such that $Q$
lies on $L_p$. Then, $\{L_p\}_{p\in \K_+^Q}$ are the
exceptional components going through $Q$. From
\ref{Cartier}, we get the following local criterion.

\begin{cor}
\label{localcriterion} Let $Q$ be a point in the
exceptional locus of $X$. If $C\subset S$, denote by $C_Q$
the curve on $S$ composed of the branches $\gamma$ of $C$
whose strict transform $\widetilde{\gamma}$ on $X$ goes
through $Q$. Then, $\widetilde{C}$ is locally principal in
a neighbourhood of $Q$ if and only if
\[\LL_{C_Q} \in \bigoplus_{u\in \K_+^{Q}} \mathbb{Z} \mathcal{L}_u.\]
\end{cor}

\begin{proof}
From the definition of $C_Q$ the germs of $\widetilde{C}$
and $\widetilde{C_Q}$ at $Q$ are equal. Moreover, $Q$ is
the unique point in the intersection of $\widetilde{C}_Q$
with the exceptional locus of $X$. Therefore,
$\widetilde{C}$ is locally principal near $Q$ if and only
if $\widetilde{C_Q}$ is Cartier. By \ref{Cartier}, this is
the case if and only if $\LL_{C_Q}\in \bigoplus_{p\in
K\_+}\mathbb{Z}\mathcal{L}_p$. By \ref{intersect/strict},
$|\widetilde{C_Q}\cdot L_p|_X=[\widetilde{C_Q},\cdot
L_p]_Q>0$ if $p\in \K_+^Q$ and zero, otherwise. The claim
follows.
\end{proof}

From \ref{Cartier} and \ref{localcriterion}, an algorithm
providing a test to verify if the strict transform on $X$
of a given curve is Cartier or locally principal at some
point is deduced. If $C\subset S$ is a curve, then the
cluster ${\mathcal Q}_{\sC}$ (defined in \ref{Cartier}) is
consistent and equivalent to the (weighted) cluster $\T_C$
defined by the system of virtual values given by
\[v_p=\left \{ \begin{array} {ll} v_p(C) & \mbox{ if }p \in
\K_+ \\ 0 & \mbox{ otherwise.}
\end{array} \right. \]
In order to know the dicritical points of ${\mathcal
Q}_{\sC}$, it is enough to unload $\T_C$. By \ref{Cartier},
$\widetilde{C}$ is Cartier if and only if $({\mathcal
Q}_{\sC})_+\subset \K_+$. To study the local principality
of $\widetilde{C}$ at some singularity $Q$, proceeding
analogously with $\Q_{\sC_Q}$ is enough. Forthcoming
\ref{Ex:1} and \ref{Ex:2} provide examples of this.

\begin{rem}
A curve on $X$ containing no exceptional component is the
strict transform of some curve on $S$, so \ref{Cartier}
actually gives criteria for the principality of any such
curves. The principality of exceptional curves can be
studied componentwise in terms of curves with no
exceptional component, taking into account \ref{mumf}: for
and $p\in \K_+$, write $\gamma_p$ for a generic curve going
sharply through $\K(p)$ and missing all points in $K$ after
$p$. Then, the principality of a curve
\[\widetilde{C}+\LL=\widetilde{C}+\sum_{p\in \K_+}a_p L_p\] is equivalent to the
principality of the strict transform on $X$ of $C+C_{\LL}$,
where
\[C_{\LL}=\sum_{p\in \K_+} a_p C^p,\]
and $C^{p}=\sum_{E_p\cap E_q\neq \emptyset} \gamma_q$.
Indeed, the curves $C^p$ are defined in this way so that
the equality $D_{\widetilde{C^p}}=D_{L_p}$ holds, and we
can apply \ref{mumf}. Now, if
$\widetilde{C}+\widetilde{C_{\LL}}$ is Cartier or not can
be decided by using the equivalent assertions of
\ref{Cartier}.
\end{rem}

\section{Some consequences and examples}

Given $C$ on $S$, let $\gamma_1,\ldots, \gamma_s$ be the
branches of $C$ at $O$ and for each $i$ denote by $p_i$ the
first non-singular point on $\gamma_i$ and not in $K$.
Write $K'$ for  the minimal cluster containing $K$ and the
points $p_1,\ldots, p_s$. Since the sheaf $I{\mathcal
O}_{S_{K'}}$ is invertible, there is a morphism
$g:S_{K'}\longrightarrow X$ induced by the universal
property of the blowing-up. Moreover, $g$ factors through
$f:S_K\rightarrow X$ because this is the minimal resolution
of $X$. Since $K'$ contains $p_1,\ldots, p_s$, the strict
transform of $C$ on $S_{K'}$ is non-singular, while the
total transform has normal crossings only:
\[\xymatrix{  {S_{K'}} \ar[r] \ar @/^2pc/[rr]^{g} \ar[rrd]_{\pi_{K'}} & S_K \ar[r]^{f} \ar[rd]^{\pi_K}
& {X} \ar[d]^{\pi_I} \\ & & {S} }\]

\noindent \emph{Remark on the notation.} From now on, we
are working on the surface $S_{K'}$ rather that on $S_K$.
For the sake of simplicity in the notation, we will make a
slight abuse of language and write $S'$ for $S_{K'}$,
$\widetilde{C}'$ and $D_{\widetilde{C}}$ for the strict
transform and the exceptional part of the total transform
of $\widetilde{C}$ on $S'$, and also $\E_J$ and $\E_C$ for
the exceptional parts of $J\in \II$ and $C$ on $S'$,
respectively. Recall from \ref{Cartier} that $J_C=\{g\in R
\mid v_p(g)\geq v_p(C), \forall p\in \K_+\}$ and
$\Q_{\sC}=BP(J_C)$.

\vspace{3mm} During the proof of \ref{Cartier} we have
proved the following fact relating the coefficients of
$\LL_C$ in the base $\{\LL_p\}_{p\in \K_+}$ to the
intersection product of $\widetilde{C}$ with $\{L_p\}_{p\in
\K_+}$.

\begin{thm} \textbf{\emph{(corollary of the proof of \ref{Cartier})}}
\label{intersect/strict} If $\LL_C=\sum_{u\in \K_+}a_u
\mathcal{L}_u$, then  \[a_u=|\widetilde{C}\cdot L_u|_X\geq
0.\]Moreover, if $\widetilde{C}$ is Cartier, then
$J_C=\prod_{u\in \K_+}I_u^{a_u}$ is the (Zariski)
factorization of $J_C$. In particular, $J_C{\mathcal
O}_X={\mathcal O}_X(-\LL_C)$.
\end{thm}

\begin{cor}
\label{hypersharplyCartier} Let $I=\prod_{p\in
\K_+}I_p^{\alpha_p}$ be the (Zariski) factorization of $I$
and let $Q$ be a point in the exceptional locus of $X$. If
$C$ goes sharply through $\K_Q$, then $\widetilde{C}$ is
Cartier on $X$ and $|\widetilde{C}\cdot L_p|_X=\alpha_p$,
for all $p\in \K_+$.
\end{cor}

\begin{proof}
We already know that for such a curve $C$, $\LL_C=\LL_I$.
Then, it is enough to apply \ref{Cartier} and
\ref{intersect/strict}.
\end{proof}

The following proposition is technical and will be useful later on.

\begin{prop}\label{hil}
\begin{itemize}
\item[(a)] If $\widetilde{C}$ is Cartier, then $D_{\widetilde{C}}=\E_C-\E_{{\mathcal Q}_{\sC}}$.

\item[(b)] If $C_1,C_2$ are curves on $S$, then
\[|\widetilde{C_1}\cdot \widetilde{C_2}|_X=[C_1,C_2]_O-[\Q_{\sC_1},C_2]_O.\]
\end{itemize}
\end{prop}

\begin{proof}
Since $\pi_I^*(C)=\widetilde{C}+\LL_C$ and
$\pi_{K'}^*(C)=\widetilde{C}+\E_C$, we infer that
\begin{eqnarray}
\label{auxDC} \E_C & = & D_{\widetilde{C}}+g^*(\LL_C).
\end{eqnarray}
Let $B$ be a curve going sharply through ${\mathcal
Q}_{\sC}$ and missing the points in $K'\setminus K$. Then,
$\LL_B=\LL_C$ and $\widetilde{B}$ goes through no
singularities of $X$ and shares no points with
$\widetilde{C}$. Thus, $\pi_{K'}^*(B)=
\widetilde{B}+g^*(\LL_C)$ and so, $\E_{{\mathcal Q}_{\sC}}
= g^*(\LL_C)$. (a) follows from (\ref{auxDC}) above.

Now, by the projection formula applied to $\pi_{K}$, we
have
\begin{eqnarray}
\label{astaux} |\widetilde{C_1}\cdot \widetilde{C_2}|_X =
|(\widetilde{C_1}+\E_{C_1})\cdot \widetilde{C_2}|_{S'}-
|\E_{{\mathcal Q}_{\sC}}\cdot \widetilde{C_2}|_{S'}
\end{eqnarray}
As above, let $B_1$ be a curve going sharply through ${\mathcal Q}_{\sC_1}$ and such that
$\widetilde{B_1}$ shares no point with $\widetilde{C_2}$ on $S'$. Then,
\begin{eqnarray*}
[{\mathcal Q}_{\sC_1},C_2]_O =
[B_1,C_2]_O=|(\widetilde{B_1}+\E_{{\mathcal
Q}_{\sC_1}})\cdot \widetilde{C_2}|_{S'} = |\E_{{\mathcal
Q}_{\sC_1}}\cdot \widetilde{C_2}|_{S'}
\end{eqnarray*} and  $[C_1,C_2]_O =
|(\widetilde{C_1}+\E_{C_1})\cdot \widetilde{C_2}|_{S'}$.
The claim of (b) is derived from (\ref{astaux}).
\end{proof}

\vspace{3mm} To close this section, we illustrate the
results in it and in the previous one with some examples.

%
%
%

\vspace{3mm}
\begin{example}\label{Ex:1}
Take $I\in \II$ with base points as on (a) of figure
\ref{Fig:1}. The dicritical points of $\K=BP(I)$ are
$p_2,p_4$ and $p_8$ and so, the surface $X=Bl_I(S)$ has
three exceptional components $L_{p_2}, L_{p_4}$ and
$L_{p_8}$. $X$ has two singularities: $Q_1$ in the
intersection of $L_{p_2}$ and $L_{p_4}$, with
$T_{Q_1}=\{p_3\}$ and $Q_2$ in the intersection of
$L_{p_2}$ and $L_{p_8}$, with
$T_{Q_2}=\{p_1,p_5,p_6,p_7\}$. We have $\mathcal{L}_{p_2} =
2L_{p_2}+2L_{p_4}+2L_{p_8}$, $\mathcal{L}_{p_4} =
2L_{p_2}+4L_{p_4}+2L_{p_8}$ and $\mathcal{L}_{p_8} =
2L_{p_4}+2L_{p_4}+9L_{p_8}$. If $C$ is a curve on $S$ with
singular points as represented on (b) of figure
\ref{Fig:1}, then $\LL_C = 6L_{p_2}+8L_{p_4}+11L_{p_8} =
\frac{7}{6}\mathcal{L}_{p_2}+{\mathcal
L}_{p_4}+\frac{5}{6}\mathcal{L}_{p_8}$. By \ref{Cartier},
$\widetilde{C}$ is not a Cartier divisor on $X$. However,
$\LL_{C_{Q_1}}=4L_{p_2}+6L_{p_4}+4L_{p_8} =
\mathcal{L}_{p_2}+\mathcal{L}_{p_4}$, so $\widetilde{C}$ is
principal in a neighborhood of $Q_1$, but not near $Q_2$
(see \ref{localcriterion}).
\begin{figure}
\begin{center}
\psfrag{A}{$Q_1$} \psfrag{B}{$Q_2$} \psfrag{X}{$X$} \psfrag{Q}{$L_{p_2}$}
\psfrag{D}{$L_{p_4}$} \psfrag{H}{$L_{p_8}$} \psfrag{Z}{$\widetilde{C}$}
\includegraphics[scale=0.6]{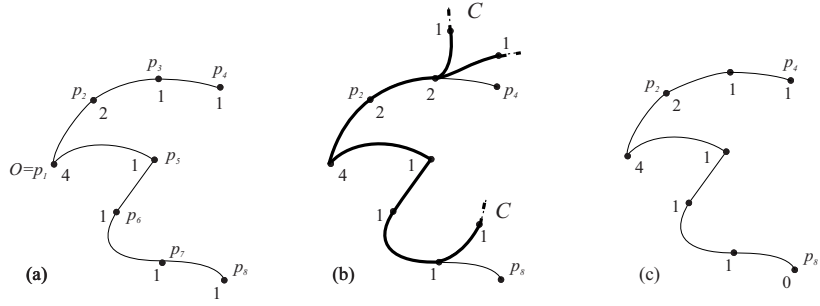}
\end{center}
\caption{\label{Fig:1}  (a) the Enriques diagram of $\K$;
(b) the singular points of $C$ in bold; (c) the Enriques
diagram of the cluster $\Q_C$ (Example \ref{Ex:1}).}
\end{figure}
\end{example}

\begin{example}
\label{Ex:2} Let $I\in \II$ with base points as on (a) of
figure \ref{Fig:2}. The dicritical points of $\K=BP(I)$ are
$p_1,p_4,p_8$ and $p_{10}$ and so, the surface $X=Bl_I(S)$
has exceptional components $L_{p_1}, L_{p_4},L_{p_8}$ and
$L_{p_{10}}$. There is only one singularity on $X$, say
$Q$, and $\mathcal{L}_{p_1} =
L_{p_1}+L_{p_4}+2L_{p_8}+2L_{p_{10}}, \mathcal{L}_{p_4} =
L_{p_1}+4L_{p_4}+4L_{p_8}+4L_{p_{10}}, \mathcal{L}_{p_8} =
2L_{p_1}+4L_{p_4}+12L_{p_8}+10L_{p_{10}},
\mathcal{L}_{p_{10}} =
2L_{p_1}+4L_{p_4}+10L_{p_8}+12L_{p_{10}}$.
\begin{figure}
\begin{center}
\includegraphics[scale=0.65]{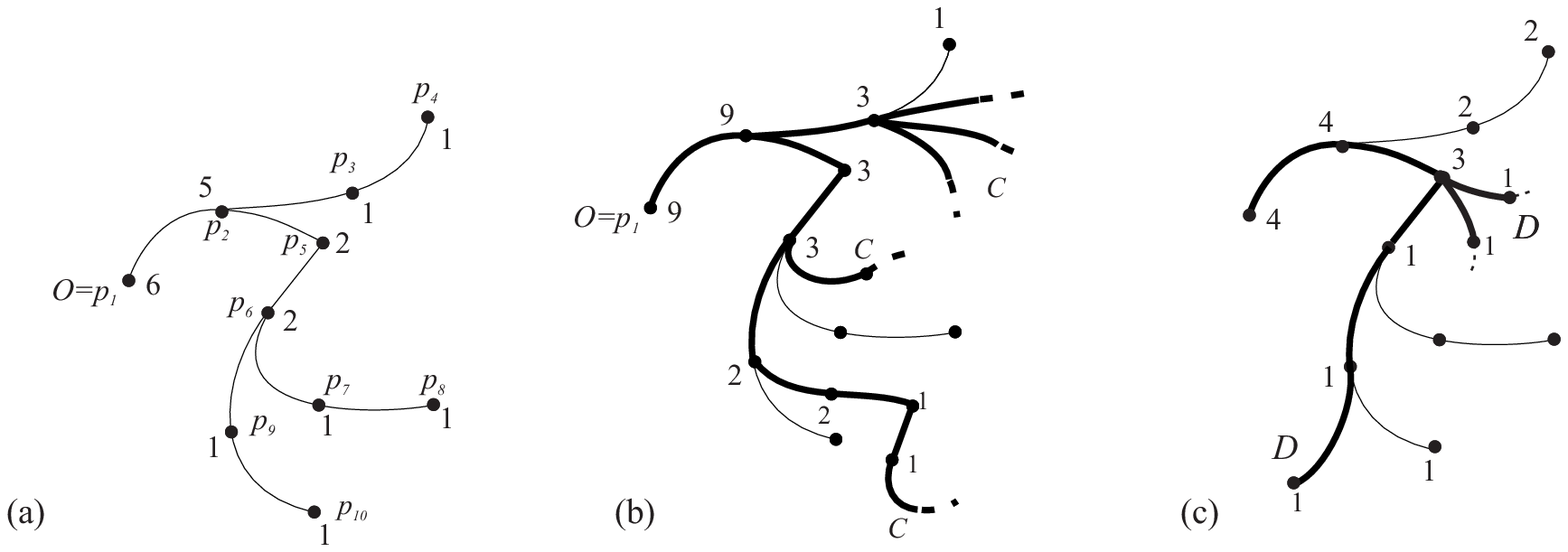}
\end{center}
\caption{\label{Fig:2}  (a) The Enriques diagrams of
${\mathcal K}$; (b) the singular points of the curves $C$;
(c) the singular points of $D$ (Example \ref{Ex:2}).}
\end{figure}
Take the curves $C$ and $D$ having singularities as (c) and
(d) of figure \ref{Fig:2}. Direct computation shows that
$\LL_C = 9 L_{p_1}+ 21 L_{p_4}+42 L_{p_8}+ 44 L_{p_{10}} =
\mathcal{L}_{p_1}+ 2 \mathcal{L}_{p_4}+ \mathcal{L}_{p_8}+
2 \mathcal{L}_{p_{10}}$, so $\widetilde{C}$ is a Cartier
divisor on $X$ and $|\widetilde{C}\cdot L_{p_1}|_X=1$,
$|\widetilde{C}\cdot L_{p_4}|_X=2$, $|\widetilde{C}\cdot
L_{p_8}|_X=1$ and $|\widetilde{C}\cdot L_{p_{12}}|_X=2$.
The cluster ${\mathcal Q}_{\sC}$ is represented on (b) of
figure \ref{Fig:2}. We have that $[C,D]_O=88$ and
$[{\mathcal Q}_{\sC},D]_O=82$. Thus, in virtue of (b) of
\ref{hil}, $|\widetilde{C}\cdot \widetilde{D}|_X=6$.
\end{example}

\vspace{2mm}
\section{On Cartier divisors containing a given Weil divisor on $X$}

In this section, we are interested in studying the
(effective)  Cartier divisors containing a given curve on
$X$. For the sake of simplicity in the exposition, we
restrict our study to curves containing no exceptional
components on $X$ and whose strict transform on $S'$ is
already non-singular. Moreover, we add a minimality
condition concerning  the values of these curves relative
to the divisorial valuations $\{v_p\}_{p\in T_Q^C}$. To
this aim, recall that if $A$  is a curve on $X$, $A$ is a
Cartier divisor if and only if $D_{A}$ is a divisor on $S'$
(see \ref{mumf}). When $A$ is not Cartier, it is also
possible to associate a divisor $\overline{D}_{A}$ on $S'$
to $A$ such that
\[ |E_p\cdot \overline{D}_A|_{S'} \leq |E_p\cdot D_{A}|_{S'} \quad
\mbox{for any exceptional component $E_p$ for $g$}\]and
minimal with this property. This divisor $\overline{D}_A$
can be easily computed by  Laufer's method; moreover, $A$
is Cartier if and only if $D_{A}=\overline{D}_{A}$ (see
\cite{Giraud} \S 1 for details).

\vspace{3mm} Fix $C$ be a curve on $S$. An effective
Cartier divisor $C'$ containing $\widetilde{C}$ will be
called \emph{$v$-minimal} (relative to $C$) if
\begin{itemize}
\item[(i)] $D_{C'}=\overline{D}_{\widetilde{C}}$;
\item[(ii)] the strict transform of $C'$ on $S'$ is non-singular and $g^*(C')$ has normal crossings only;
\item[(iii)] $C'$ contains no exceptional components on $X$.
\end{itemize}
By using partial unloading (see Appendix), we  are going to
construct flags of clusters that will allow us to give a
complete description of the $v$-minimal Cartier divisors
containing $\widetilde{C}$. These flags will play a basic
role in the forthcoming sections since, as it will be shown
there, they carry deep information relative to
$\widetilde{C}$.

Recall that $K'$ is the minimal cluster containing $K$ and
the first non-singular points $p_i$ of the branches of $C$
which are not in $K$. Fix non-negative integers
$\mathbf{m}=\{m_p\}_{p\in \K_+}$, and write
$\K^{\mathbf{m}}$ for the cluster whose underlying set of
points is $K'$ and whose excess at $p$ is $m_p$ if $p\in
\K_+$ and $0$, otherwise. Take the sequence of clusters
\begin{eqnarray} \label{cahin}
\T_0=\K^{\mathbf{m}} \prec \T_1 \prec \ldots \prec \T_j
\prec \ldots \hspace{4mm} \mbox{ with }\T_j=(K',\tau^{j})
\end{eqnarray}
defined as follows: for $j\geq 0$, as far as there exists
some $p_i$ such that $\tau^{j}_{p_i}=0$, $\T_{j+1}$ is the
cluster obtained from $\T_j$ by increasing its multiplicity
at $p_i$ by one and performing partial unloading relative
to the set $\K_+$. Note that as long as there exists such a
$p_i$, $H_{\T_j}\supsetneq H_{\T_{j+1}}$, and moreover
\begin{eqnarray}\label{SH}
dim_{\mathbb{C}} (\dfrac{H_{\T_j}}{H_{\T_{j+1}}})=1.
\end{eqnarray}
A flag of clusters as in  (\ref{cahin}) will be called a
\emph{flag for} ($\K$,$C$) (an \emph{$\mathbf{m}$-flag} if
we want to precise the original excesses of $\T_0$). For
simplicity in the notation and if no confussion may arise,
we will write $\{T_j\}$ for the clusters appearing in such
a flag.

\begin{lem}\label{finsteps}
After finitely many steps, this procedure stops giving rise
to a cluster $\T_n$ such that
\begin{itemize}
\item[1.] $\tau^n_{p_i}=1$, for $i=1,\ldots,s$;
\item[2.] $\LL_{\T_j}=\sum_{p\in \K_+}m_p\LL_p$, for $j=0,\ldots,n$.
\end{itemize}
and non-negative integers
$\boldsymbol{\omega}=\{\omega_p\}_{p\in \K_+}$ such that
$\rho^{\T_n}_p=m_p-\omega_p$ for $p\in \K_+$. Moreover, the
number $n$ and the integers $\boldsymbol{\omega}$ depend
only on $\K$ and the points $p_1,\ldots,p_s$, and not on
$\mathbf{m}$. \vspace{3mm}
\end{lem}

\begin{proof}
Write
$\overline{\K^{\mathbf{m}}}=\K^{\mathbf{m}}+\sum_{i=1}^m
\K(p_i)$. First of all, we show by using induction on $j$
that
\begin{eqnarray}\label{for:inclusiondelta}
H_{\overline{\K^{\mathbf{m}}}}\subset H_{\T_j}
\end{eqnarray}
for every $\T_j$ defined as above. This is clear for $j=0$.
Assume it is also true for some $j\geq 0$. To prove that
$H_{\overline{\K^{\mathbf{m}}}}\subset H_{\T_{j+1}}$, it is
enough to show that the virtual transform relative to the
multiplicities of  $\widetilde{\T_j}$ of a generic curve
going through $\overline{\K^{\mathbf{m}}}$ goes through the
point $p_i$. But this is clear since any curve going
sharply through $\overline{\K^{\mathbf{m}}}$ goes
\emph{effectively} through $p_i$. Because of (\ref{SH}),
this procedure stops after finitely many steps. Hence, we
obtain a (not necessarily consistent) cluster $\T_n$ such
that
\begin{eqnarray*}\label{aha}
\tau^n_{p_i}\geq 1 \ \mbox{ for }i=1,\ldots,s.
\end{eqnarray*}
Notice by the way that $n$ is the codimension of $\T_n$ in
$\T_0$.

Now, assume that we have a couple of (not necessarily
consistent) clusters $\T^{(1)}=(K',\tau^{(1)})$ and
$\T^{(2)}=(K',\tau^{(2)})$ such that $H_{\T^{(i)}}\subset
H_{\K^{\mathbf{m}}}$ and $\tau^{n}_{p_i}\geq 1$ for
$i=1,\ldots,s$. Define a new cluster
$\T^{(0)}=(K',\tau^{(0)})$ as follows:
take\[\textbf{v}_{(0)}=(v^0_p)_{p\in K}\] where
\[v^0_p=\min\{v_p^{\T^{(1)}},v_p^{\T^{(2)}}\},\]and define
virtual multiplicities for $K'$ by
\[\tau^{(0)}_p=\left \{ \begin{array}{ll}\textbf{1}_p^t\textbf{P}_K\textbf{v}_{(0)}
& \mbox{ if }p\in K \\ \overline{\tau}_p & \mbox{ if }p\in
K'\setminus K.
\end{array} \right.\]where  $\{\overline{\tau}_p\}_{p\in K'}$ are the virtual
multiplicities of  $\overline{\K^{\mathbf{m}}}$. Clearly,
$H_{\T^{(1)}},H_{\T^{(2)}}\subset H_{\T^{(0)}}$ and
$\tau^{(0)}_p=\overline{\tau}_p$ for $p\in K'\setminus K$.
Since
$v_p^{\T^{(0)}}=\min\{v_p^{\T^{(1)}},v_p^{\T^{(2)}}\}\geq
v_p^{\K^{\mathbf{m}}}$ for any $p\in K$, we have
$H_{\T^{(0)}} \subset H_{\K^{\mathbf{m}}}$. By using
Artin's trick as in the proof of \ref{partunload}, we see
that $\rho^{\T^{(0)}}_p\geq 0$ for $p\in K'\setminus \K_+$.
This proves the uniqueness of a minimal cluster
$\T'=(T,\tau')$ with $H_{\T'}\subset H_{\K^{\mathbf{m}}}$
and $\tau^n_{p_i}\geq 1$ for $i=1,\ldots,s$, and shows
actually that $\tau'_{p_i}=1$ for $i=1,\ldots,s$. From the
way the cluster $\T_n$ has been constructed, necessarily
$\T'=\T_n$, and so $\T_n$ does not depend on the choices
done when constructing the flag (\ref{cahin}). In virtue of
\ref{tech}, if $j\geq 1$ and $p\in T\setminus \K_+$,
$\tau^{j-1}_p$ does not depend on the excesses
$\{\rho^{\T_j}_q\}_{q\in \K_+}$. By  induction, the
multiplicities $\{\tau^{j}_p\}_{p\notin \K_+}$ do not
depend on $\mathbf{m}$, and neither the difference
$\omega_p^j=\rho_p^{\T_{j-1}}-\rho_p^{\T_{j}}$. Hence, the
numbers $\omega_p=\sum_{j=1}^n\omega^{j}_p$ ($p\in \K_+$)
are independent of $\mathfrak{m}$ and
$\rho^{\T_n}_p=m_p-\omega_p$ as claimed.
\end{proof}

\begin{rem}\label{turner}
Notice that $\mathbf{m}\geq \boldsymbol{\omega}$ (this
meaning that $m_p\geq\omega_p, \forall p\in \K_+$) if and
only if all the clusters $\T_0,\T_1,\ldots,\T_n$ are
consistent. In this case, partial unloading relative to
$\K_+$ equals usual unloading.
\end{rem}

Write $\Q'_{\sC}=(K',e_C)$ for the cluster obtained by
taking  as virtual multiplicities the effective
multiplicities of $C$. If $\widetilde{C}$ is already
Cartier, the integers $\{\omega_p\}_{p\in \K_+}$ given by
\ref{finsteps} equal the effective values $\{v_p(C)\}_{p\in
\K_+}$. In fact, in this case we have that
$\Q_{\sC}=\T_0^{\boldsymbol{\omega}}$ and
$\Q'_{\sC}=\T_n^{\boldsymbol{\omega}}$.

\vspace{2mm} We state the following easy lemma for future
reference.

\begin{lem}
\label{laia} Let $\mathbf{m}\geq \boldsymbol{\omega}$ and
$\{\T_i\}_{i=0,\ldots,n}$ an $\mathbf{m}$-flag for
($\K$,$C$). For any curve $B$ on $S$, we have that
\begin{itemize}
\item[(a)] the value of $[\T_j,B]_O-[\T_0,B]_O$ does not depend on $\mathbf{m}$.
\item[(b)] if $\widetilde{B}$ shares no points on $X$ with $\widetilde{C}$, $[\T_j,B]_O=[\T_0,B]_O$ fopr any $j\geq 0$. In particular,
\[\sum_p\tau^{j}_p e_p(B)=\sum_p\tau^{0}_p e_p(B).\]
\end{itemize}
\end{lem}

\begin{proof}
If $C_j$ is a curve going sharply through $\T_j$, we have
that $C_j^*=\widetilde{C_j}+\LL_{\mathbf{m}}$. Then, by the
projection formula applied to $\pi$,
$[\T_j,B]_O=|(\widetilde{C_j}+\LL_{\mathbf{m}})\cdot
\widetilde{B}|_X$. The first claim follows immediately.
Now, notice that for a generic  curve going through $\T_0$,
$|\widetilde{C_0}\cdot \widetilde{B}|_X=0$. If
$\widetilde{B}$ shares no points on $X$ with
$\widetilde{C}$ and $j>0$, then $|\widetilde{C_j}\cdot
\widetilde{B}|_X=0$ also. The second assertion is derived
now from the first one and (b) of \ref{hil}.
\end{proof}

The following proposition is the main result of this
section and gives a complete description of the $v$-minimal
Cartier divisors containing $\widetilde{C}$.

\begin{prop} \label{minCardiv}
Take $\boldsymbol{\omega}$ as in \ref{finsteps}.
\begin{itemize}
\item[(a)] If $\mathbf{m}\geq \boldsymbol{\omega}$, the system of virtual multiplicities $\{\tau^{n}_p-e_p(C)\}_{p\in K'}$ is consistent.

\item[(b)]Every $v$-minimal Cartier divisor for $\widetilde{C}$
has the form $\widetilde{C}+\widetilde{B}$, where $B$ is a
generic curve going through ${\mathcal
C}_{\mathbf{m}}=(K',\tau-e_C)$, for some $\mathbf{m}\geq
\boldsymbol{\omega}$.

\item[(c)] $\overline{D}_{\widetilde{C}}=
{\E}_{\T_n}-{\E}_{\T_0}=\sum_{p\in K'}n_pE_p$, where $n_p$
is the number of unloadings performed on each $p\in K'$.

\end{itemize}
\end{prop}

\begin{rem}\label{CartC2}
After \ref{minCardiv}, one should think of the strict
transform of the generic curves going through ${\mathcal
C}_{\mathbf{m}}$ (for some $\mathbf{m}\geq
\boldsymbol{\omega}$) as the curves to be added to
$\widetilde{C}$ to obtain ($v$-minimal) Cartier divisors on
$X$.
\end{rem}

\begin{proof}
Since $\mathbf{m}\geq \boldsymbol{\omega}$, the cluster
$\T^{\mathbf{m}}_n$ is consistent (see \ref{turner}), and
from its construction, the points $p_1,\ldots,p_s$ have
virtual multiplicity one at it and are maximal in $K'$. It
follows that $\rho^{\T^{\mathbf{m}}_n}_{p_i}=1$ for $i\in
\{1,\ldots,s\}$. On the other hand, from its own definition
the cluster $\Q'_{\sC}$ has excess equal to one at
$\{p_i\}_{i=1,\ldots,s}$,  and zero at the remaining
points. Since the excesses of $\mathcal{C}_{\mathbf{m}}$
are the difference between the excesses of
$\T^{\mathbf{m}}_n$ and $\Q'_{\sC}$, we infer that
${\mathcal{C}}_{\mathbf{m}}$ is consistent.

Let $A$ be an effective Cartier divisor on $X$ containing
$\widetilde{C}$ and without exceptional components. There
exists some curve $B$ on $S$ such that $A$ is the strict
transform of $C_B=C+B$. By \ref{Cartier}, we have that
\begin{eqnarray}\label{LauxL2}
\LL_{C_B}=\sum_{p\in \K_+}m_p\mathcal{L}_p \quad \mbox{with
} \{m_p\}\subset \mathbb{Z}_{\geq 0}.
\end{eqnarray}
Write $\mathbf{m}=\{m_p\}_{p\in \K_+}$. By \ref{turner}, we
have that $\mathbf{m}\geq \boldsymbol{\omega}$ and so,  the
clusters in any $\mathbf{m}$-flag \[\T_0^{\mathbf{m}}\prec
\ldots \prec \T_n^{\mathbf{m}}\]are consistent. Now, $C_B$
goes through $\T_0^{\mathbf{m}}$ and $e_{p_i}(C_B)=1$ for
$i=1,\ldots,s$. By the minimality of $\T_N^{\mathbf{m}}$
(see \ref{finsteps}), it follows that $C_B$ goes through
it. Thus, $C_B$ goes virtually through  $\T_n^{\mathbf{m}}$
and $\Q_{C_B}=\T_0^{\mathbf{m}}$. Hence, by (a) of
\ref{hil},
\begin{eqnarray*}
D_{C_B} = \E_{C_B}-\E_{\T_0^{\mathbf{m}}} \geq
\E_{\T_n^{\mathbf{m}}}- \E_{\T_0^{\mathbf{m}}},
\end{eqnarray*}and the equality holds if and only if $C_B$ goes
through $\T_n^{\mathbf{m}}$ with effective multiplicities
equal to the virtual ones. This is equivalent to say that
the curve $B$ goes through ${\mathcal C}_{\mathbf{m}}$ with
effective multiplicities equal to the virtual ones. By
definition of the cluster $K'$, $\widetilde{C_B}'$
intersects transversally the exceptional divisor of $f_C$.
Therefore, the strict transform $\widetilde{A'}^{S_{K_C}}$
intersects transversally the exceptional divisor of $f_C$
if and only if for each $p\in K'\setminus \K_+$, $C_B$ has
$\rho^{\T_n^{\mathbf{m}}}_p$ branches through $p$ and
misses all points after $p$ in $K'$. Therefore, we see that
$A$ is a $v$-minimal Cartier divisor containing
$\widetilde{C}$ if and only if $B$ goes through ${\mathcal
C}_{\mathbf{m}}$ with effective multiplicities equal to the
virtual ones and has $\rho^{\T_n^{\mathbf{m}}}_p$ branches
through $p$ and misses all points after $p$ in $K$. This
completes the proof of (b).

Moreover, if $A$ is $v$-minimal, then
$D_{A}=\E_{\T_n^{\mathbf{m}}}- \E_{\T_0^{\mathbf{m}}}$ and
in particular, $\overline{D}_{\widetilde{C}}=
\E_{\T^{\mathbf{m}}_n}-\E_{\T^{\mathbf{m}}_0}$. Now, from
the definition of $\T_n^{\mathbf{m}}$, we know that
$\tau^{\mathbf{m}}_p=e_p(C)$ for every $p\in K'\setminus
K$. Therefore, $p\notin K$ is a dicritical point of
$\T_n^{\mathbf{m}}$ if and only if $p=p_i$ for some $i$.
From its own definition, it follows that
$\mathcal{C}_{\mathbf{m}}$ has no dicritical points out of
$K$. This gives (c) and completes the proof of the
proposition.
\end{proof}

We finish with an example

\begin{example}\label{Ex:4}
Take $I\in \II$ as in Example \ref{Ex:2}. On (a) of figure
\ref{Fig:6}, the clusters $\T_0$ and $\T_n$ of a
$\mathbf{m}$-flag for ($\K$, $C$) are represented. The
excesses of $\T_n$ at $p_2,p_4$ and $p_8$ are
$\rho^{\T_n}_{p_2}=m_{p_2}-2$,
$\rho^{\T_n}_{p_4}=m_{p_4}-1$,
$\rho^{\T_n}_{p_8}=m_{p_8}-2$, and so, $\omega_{p_2}=2$,
$\omega_{p_4}=1$ and $\omega_{p_8}=2$. By taking
$\mathbf{m}=\boldsymbol{\omega}$, the virtual values of
$\T_n$ and $\T_0$ are
$\mathbf{v}_{\T_n}=\{8,10,12,12,10,20,22,22\}$ and
$\mathbf{v}_{\T_0}=\{7,10,11,12,9,18,20,22\}$. By
\ref{minCardiv}, the $v$-minimal Cartier divisors on $X$
containing $\widetilde{C}$ are the curves
$\widetilde{C}+\widetilde{B}$, where $B$ is a generic curve
going through ${\mathcal C}_{\mathbf{m}}$. By (b) of
\ref{minCardiv}, $\overline{D}_{\widetilde{C}} =
\E_{\T_n}-\E_{\T_0}
 = E_{p_3}+E_{p_1}+E_{p_5}+2E_{p_6}+E_{p_7}$, so
$\overline{D}_{\widetilde{C}}=D^{Q_1}_{\widetilde{C}}+
\overline{D}^{Q_2}_{\widetilde{C}}$.

\begin{figure}
\begin{center}
\includegraphics[scale=0.75]{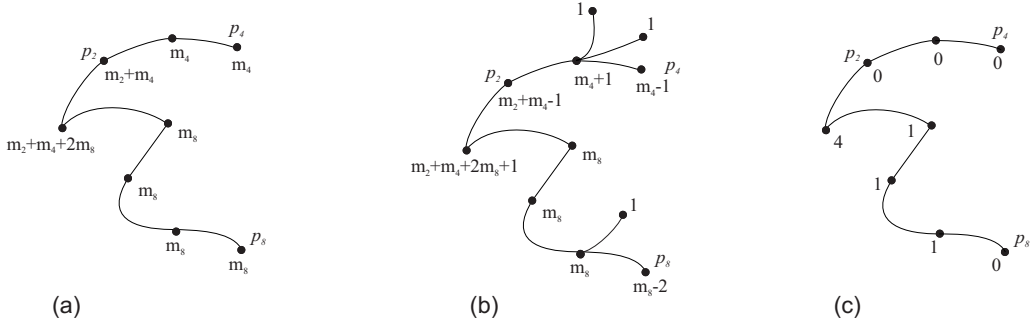}
\end{center}
\caption{\label{Fig:6} (a) the Enriques diagram of
$\T_0^{\mathbf{m}}$; (b) the Enriques diagram of
$\T_n^{\mathbf{m}}$; (c) the Enriques diagram of ${\mathcal
C}_{\boldsymbol{\omega}}$ (Example \ref{Ex:4}).}
\end{figure}
\end{example}

\vspace{2mm}
\section[On the order of singularity of curves]{The order of singularity of curves on $X$}
For a given curve $C$ on $S$, write
$\delta_X(\widetilde{C})=\sum_{Q}\delta_Q(\widetilde{C})$
and call this value the \emph{total order of singularity of
$\widetilde{C}$}. Our aim in this section is to get some
formulas for the order of singularity of curves on $X$
without exceptional support. Then, we will use them to
compute the semigroup of a branch going through a
sandwiched singularity

\begin{prop}
\label{delta/general} Let $C$ be a curve on $S$ and take a
flag $\{\T_i\}_{i=0,\ldots,n}$
for $(\K,C)$. Then,
\[\delta_X(\widetilde{C})=[\T_n,C]_O-[\T_0,C]_O-
n.\]
\end{prop}

In virtue of \ref{laia} the right term of the above
equality does not depend on the particular flag
$\{\T_i\}_{i=0,\ldots,n}$ chosen.

\vspace{3mm} To prove this result we need a technical
lemma.

\begin{lem}
\label{div/canonical} Let $\T=(T,\tau)$ be a a consistent
cluster, and write $\mathbb{K}_{T}$ for the canonical
divisor on the surface $S_T$ obtained by blowing up all the
points of $T$. If $A$ is a curve on $S$ and $\E_A$ is the
exceptional part of its total transform on $S_T$, then
\[|\E_{A}\cdot \mathbb{K}_{T}|_{S_T}=-\sum_{p\in T}e_p(A).\]
\end{lem}

\begin{proof}
The adjunction formula (see \cite{Hartshorne} V.1.5) says
that
\[|E_p \cdot \mathbb{K}_{T}|_{S_T}=-2-E_p^2 \quad \mbox{for each $p\in T$.}\]

Write $\T'$ for the cluster obtained by taking
$\{\tau_p-1\}_{p\in T}$ as the system of virtual
multiplicities. Let $B,B'$ be curves going sharply through
$\T$ and $\T'$, respectively. Direct computation using
(\ref{for:int/comp}) and that $E_p^2=-\sharp\{q\in T \mid
q\rightarrow p\}-1$ shows that \[|E_p\cdot
\widetilde{B'}^{T}|_{S_T}- |E_p\cdot
\widetilde{B}^{T}|_{S_T}=-2-E_p^2.\]Therefore, by the
projection formula,
\[|\E_{A}^{T}\cdot \mathbb{K}_{T}|_{S_T}=
|\E_{A}^{T}\cdot \widetilde{B'}^{T}|_{S_T}-
|\E_{A}^{T}\cdot
\widetilde{B}^{T}|_{S_T}=[A,B']_O-[A,B]_O=-\sum_{p\in
T}e_p(A),\]the last equality by Noether's formula.
\end{proof}

\vspace{3mm}

\noindent \textsc{Proof of \ref{delta/general}.} Write
$\mathbb{K}$ for the canonical divisor on $S'$ For the sake
of simplicity, write $D$ and $\overline{D}$ for
$D_{\widetilde{C}}$ and $D\overline{}_{\widetilde{C}}$,
respectively. Add the subindex $Q$ to denote the
exceptional part of $D$ and $\overline{D}$ contracting to
some $Q\in X$. By Corollary 2.1.2. of \cite{Morales}, we
have the following formula for the order of singularity of
$\widetilde{C}$ at $Q$:
\begin{eqnarray}
\label{auxpredelta} \delta_Q(\widetilde{C}) & = &
\frac{1}{2}|D_Q  \cdot (\mathbb{K}-D_Q  )|_{S'}
+\frac{1}{2}|(D_Q  -\overline{D}_Q  )\cdot (D_Q
-\overline{D}_Q  -\mathbb{K})|_{S'}
= \nonumber \\
& = &  \frac{1}{2}|\overline{D}_Q  \cdot (\overline{D}_Q
-2D_Q  +\mathbb{K})|_{S'}.
\end{eqnarray}
From this and the fact that $D_Q$ and $D_{Q'}$ are disjoint
if $Q\neq Q'$, we have \[\delta_Q(\widetilde{C}) =
\frac{1}{2}|\overline{D}_Q\cdot (\overline{D}- 2D
+\mathbb{K})|_{S'}.\] Thus,
\begin{eqnarray}
\label{for:delta} \sum_{Q}\delta_Q(\widetilde{C}) & = &
\frac{1}{2}|\overline{D}\cdot
(\overline{D}+\mathbb{K})|_{S'}- |\overline{D}\cdot
D|_{S'}.
\end{eqnarray}
In virtue of (b) of \ref{minCardiv}, $\overline{D}=
\E_{\T_n}-\E_{\T_0}$. Let $B$ be a curve going sharply
through $\T_n$. Then, $\widetilde{B}$ is Cartier on $X$ and
by (a) of \ref{hil}, $D_{\widetilde{B}} = \overline{D}$. By
the projection formula applied to $f_C$,
$|\overline{D}\cdot E|_{S'}=-|\widetilde{C}'\cdot E|_{S'}$
for any exceptional divisor $E$ being contracted by
$f_{K_C}$. Therefore, we have
\begin{eqnarray}
\label{for:delta/general} &  & |\overline{D}\cdot
(\overline{D}+\mathbb{K})|_{S'}
=|(\E_{\T_n}-\E_{\T_0})\cdot
(-\widetilde{C}'+\mathbb{K})|_{S'}.
\end{eqnarray}
Now, since a generic curve through $\T_n$ (resp. through
$\T_0$) goes sharply through it and shares no points with
$B$ outside $K'$  (see Theorem 4.2.8 of \cite{CasasBook}),
we have that \[-|(\E_{\T_n}-\E_{\T_0})\cdot
B|_{S'}=[\T_0\cdot B]_O-[\T_n\cdot B]_O.\]On the other
hand, in virtue of (\ref{for:delta}),
\[|(\E_{\T_n}-\E_{\T_0})\cdot
\mathbb{K}|_{S'}=\sum_{p\in K'}(\tau^0_p-\tau^n_p).\]From
all of this, (\ref{for:delta/general}) and the Noether's
formula (see (\ref{for:Noether})), we infer that
\begin{eqnarray*}
|\overline{D}\cdot (\overline{D}+\mathbb{K})|_{S'} & = & [\T_0\cdot B]_O-[\T_n\cdot B]_O+\sum_{p\in K'}(\tau_p^{0}-\tau_p^{n}) =\\
& = & \sum_p \tau_p^{0}(\tau_p^{0}-1)-\sum_p
\tau_p^{n}(\tau_p^{n}-1),
\end{eqnarray*} the last equality because $\sum_p\tau^{0}_p\tau^{n}_p=\sum_p (\tau^{n}_p)^2$ (see (b) of \ref{hil}).
Hence,
\begin{eqnarray*}
\frac{1}{2}|\overline{D}\cdot
(\overline{D}+\mathbb{K})|_{S'} & = & \frac{1}{2}\sum_p
\tau_p^{0}(\tau_p^{0}-1)-\frac{1}{2}\sum_p
\tau_p^{n}(\tau_p^{n}-1) =
\dim_{\mathbb{C}}(\frac{H_{\T_0}}{H_{\T_n}})=n.
\end{eqnarray*}
On the other hand, since $D_{\widetilde{C'}}  =
\overline{D}$, the projection formulas for $f_C$ and
$\pi_{K_C}$ say that $|\overline{D}\cdot D|_{S'}=
-|\overline{D}\cdot \widetilde{C}'|_{S'}
=[\T_0,C]_O-[\T_n,C]_O$. The claim follows from
(\ref{for:delta}). \qed

\vspace{3mm}

From \ref{delta/general} we deduce the following formula
for the order of singularity of a Cartier divisor on $X$.

\begin{cor}
\label{delta} If $\widetilde{C}$ is a Cartier divisor on
$X$ and $d_p=e_p(C)-\tau^{0}_p$, we have that
\[\delta_{X}(\widetilde{C})=\delta_O(C)-\delta_O({\mathcal Q}_{\sC})=
\sum_{p\in K_C} \frac{d_p(d_p-1)}{2}.\]
\end{cor}

\begin{proof}
If $\widetilde{C}$ is Cartier,
$\T_n^{\boldsymbol{\omega}}=\Q'_{\sC}$ and
$\T_0^{\boldsymbol{\omega}}={\mathcal Q}_{\sC}$ (see
\ref{CartC2}). The strict transfrom of a generic curve $B$
going sharply through $\T_0$ shares no points with
$\widetilde{C}$ on $X$. Therefore, by (b) of  \ref{hil}, we
have
\[[\T_0,C]_O=[B,\T_n]_O=[B,\T_0]_O=\T_0^2\]Thus, $[\T^{\boldsymbol{\omega}}_n,C]_O=(\Q'_{\sC})^2$
and also $[\T^{\boldsymbol{\omega}}_0,C]_O={\mathcal
Q}_{\sC}^2$. Using \ref{delta/general} and the fact that
$\delta_{\T}+c(\T)=\T^2$, the first equality follows. The
second follows by direct computation since by (b) of
\ref{hil} again, $\sum_{p}d_p \tau^{0}_p=0$.
\end{proof}

Let $R_f=\frac{R}{(f)}$  be the local ring  of $C$ at $O$
and assume that $\widetilde{C}$ is Cartier. Write
$\{Q_1,\ldots,Q_n\}$ for the points of $\widetilde{C}$ in
the exceptional locus of $X$ and $J_i$ for the ideal of
$\widetilde{C}$ in ${\mathcal O}_{X,Q_i}$. Clearly, we have
a monomorphism \begin{eqnarray*} R_f & \rightarrow &
R^I_f=\prod_{i=1}^n \frac{{\mathcal O}_{X,Q_i}}{J_i}
\end{eqnarray*}induced by the monomorphisms $R\rightarrow
{\mathcal O}_{X,Q_i}$. Then, the above corollary says that
\[\dim_{\mathbb{C}}\frac{R^I_f}{R_f}=\delta_O({\mathcal
Q}_{\sC}).\] Notice that if $I=m_O$, $R^{I}_f$ is the ring
of $C$ in the first neighbourhood of $O$ introduced by
Northcott (see \cite{North55,North58}). Then, \ref{delta}
says that
$\delta_O(C)-\delta_X(\widetilde{C})=\frac{e_O(C)(e_O(C)-1)}{2}$
which is known to be the variation in the order of
singularity of $C$ after blowing up the point $O$.

\begin{example} \label{Ex:5}
Take the cluster and the curve $C$ of Example \ref{Ex:1}
(see figure \ref{Fig:1}). We have that
$[\T^{\boldsymbol{\omega}}_n,C]_O=49$ and
$[\T^{\boldsymbol{\omega}}_0,C]_O=42$. Direct computation
shows that
$\dim_{\mathbb{C}}(\frac{H_{\T^{\boldsymbol{\omega}}_0}}{H_{\T^{\boldsymbol{\omega}}_n}})=6$.
Therefore, by \ref{delta/general}, we obtain
$\delta_X(\widetilde{C})=1$. Since $Q_1$ is the only
singularity of $\widetilde{C}$, we derive that
$\delta_{Q_1}(\widetilde{C})=1$. Let $B$ be a curve going
sharply through ${\mathcal C}_{\boldsymbol{\omega}}$ and
let $C_B=C+B$. Then, as noticed above, $\widetilde{C_B}$ is
Cartier and $\Q_{\sC_B}=\T^{\boldsymbol{\omega}}_0$. By
\ref{delta},
$\delta_X(\widetilde{C_B})=\delta_O(C_B)-\delta_{\T^{\boldsymbol{\omega}}_0}
=\delta_O(C_B)-\delta_O(\T^{\boldsymbol{\omega}}_0)=33-28=5$.
In Example \ref{Ex:2} we have seen that $\widetilde{C}$ is
principal in a neighbourhood of $Q_1$, so the germs of
$\widetilde{C}$ and $\widetilde{C_B}$ in a neighbourhood of
$Q_1$ are equal, and $\delta_{Q_1}(\widetilde{C_B})=1$.
Therefore, $\delta_{Q_2}(\widetilde{C_B})=4$.
\end{example}

\vspace{2mm}
\subsection*{The semigroup of a branch on a sandwiched singularity}

Let $C$ be a branch on $S$ and let $Q$ be the point where
$\widetilde{C}$ intersects the exceptional divisor of $X$.
Here we show how the flags for $(\K,C)$ already defined
determine the semigroup of $\widetilde{C}$ on $(X,Q)$.
Recall that since $(\widetilde{C},Q)$ is analytically
irreducible, the integral closure $\overline{{\mathcal
O}_{\widetilde{C},Q}}\cong \mathbb{C}\{t\}$ is a discrete
valuation ring and the semigroup $\Sigma_Q(\widetilde{C})$
of $\widetilde{C}$ at $Q$ is
\[\Sigma_Q(\widetilde{C})=\{v_t(g) \mid g\in {\mathcal O}_{\widetilde{C},Q}\}\] where $v_t$
is the valuation corresponding to $\overline{{\mathcal
O}_{\widetilde{C},Q}}$. It is also known that
$\delta_Q(\widetilde{C})$ is the number of elements in
$\mathbb{N}\cup\{0\}\setminus \Sigma_Q(\widetilde{C})$ (see
Appendix of \cite{zartes}).

The following proposition gives a description of
$\Sigma_Q(\widetilde{C})$ in terms of differences of
intersection multiplicities of curves at $O$. \vspace{3mm}

\begin{prop} \label{semlemma}Take a flag $\T_0\prec \T_1 \ldots \prec \T_n$ for $(\K,C)$, and for $0\leq i\leq n$,
write
\[\alpha_C^{i}=[\T_i,C]_O-[\T_0,C]_O.\]Then, the
integers $\{\alpha_C^{i}\}_{0\leq i\leq n}$ are the first
$n+1$ elements of the semigroup $\Sigma_Q(\widetilde{C})$.
Moreover $\alpha_C^{n}\geq c(Q)$, the conductor of
$\Sigma_Q(\widetilde{C})$.
\end{prop}

\begin{proof}
For each $i$, let $C_i$ be a curve going sharply through
$\T_i$. Then $\LL_{C_i}\in \bigoplus_{u\in
\K_+}\mathbb{Z}\mathcal{L}_u$ and so, the strict transform
$\widetilde{C_i}$ is a Cartier divisor on $X$. Let $g_i\in
\mathfrak{m}_{X,Q}$ be a local equation for
$\widetilde{C_i}$ near $Q$, and $\overline{g_i}$ its class
in ${\mathcal O}_{\widetilde{C},Q}$. Then,
\[[\T_i,C]_O-[\T_0,C]_O=
[\widetilde{C_i},\widetilde{C}]_Q=v_t(\overline{g_i})\] is
an element in $\Sigma_Q(\widetilde{C})$. On the other hand,
$[\T_i,C]_O<[\T_{i+1},C]_O$ for each $i$. From this it
follows that $\alpha_C^{i}\in \Sigma_Q(\widetilde{C})$ and
\begin{eqnarray*}
\sharp [0,\alpha_C^{n}]\cap \Sigma_Q(\widetilde{C}) \geq
n+1.
\end{eqnarray*}
Since $\delta_Q(\widetilde{C})$ is the number of elements
in the complement of $\Sigma_Q(\widetilde{C})$,
\ref{delta/general} implies that the above inequality is
actually an equality.

Now, if $j\geq \alpha_C^{n}$, take a curve $C'$ going
sharply through $\T_n$ and sharing exactly $j
+[\T_0,C]_O-[\T_n,C]_O$ points with $C$ outside $K_C$.
Then, $\LL_{C'}=\LL_{\T_n}$ and so, $\widetilde{C'}$ is
Cartier (see \ref{Cartier}). From Noether's formula, we
infer that
\begin{eqnarray*}
[C',C]_O & = & \sum_{p\in K_C}e_p(C')e_p(C)+(j
+[\T_0,C]_O-[\T_n,C]_O) = j +[\T_0,C]_O.
\end{eqnarray*}Hence, $j=[C',C]_O-[\T_n,C]_O=
v_t(g')$ belongs to $\Sigma_Q(\widetilde{C})$ and so,
$c(\Sigma_Q(\widetilde{C}))\leq \alpha_C^n$.
\end{proof}

The above proposition provides an easy way to compute the
semigroup of $\widetilde{C}$ once a flag for $(\K,C)$ has
been computed: the differences
$\alpha_C^{i}=[\T_{i},C]_O-[\T_0,C]_O$ for $i=0,\ldots,
N-1$ provide the first elements of
$\Sigma_Q(\widetilde{C})$; then add all the integers
greater or equal than $\alpha_C^n$.

\begin{rem} The bound for the conductor of \ref{semlemma} is far from being sharp, and in general the semigroup
$\Sigma_Q(\widetilde{C})$ is not symmetric as shown in the
following example (so the curves $\widetilde{C}$ needs
 not to be Gorenstein, see \cite{kunz}
for details).
\end{rem}

\begin{example}
\label{Ex:6} Take a cluster $\K$ and a curve $C$ as shown
in figure \ref{Fig:7}. By means of the algorithm explained
in section 5, we compute the clusters $\T_n$ and $\T_0$ for
some $\mathbf{m}$-flag for $(\K,C)$. The virtual
multiplicities of $\T_n$ and $\T_0$ are respectively,
$\boldsymbol{\nu}^{\T^{\mathbf{m}}_n}=\{12,1,0,8,1,1,0,0,0,5,2,2,1,1,1\}$,
$\boldsymbol{\nu}^{\T^{\mathbf{m}}_0}=\{11,1,1,4,4,2,2,2,2\}$.
The cluster $\T_n$ is shown on (c) of figure \ref{Fig:7}.
Then, $[\T_n,C]_O=176$ and $[\T_0,C]_O=105$. Also, we have
$\dim_{\mathbb{C}}(\frac{H_{\T_0}}{H_{\T_n}})= 141-100=41$.
Thus, by \ref{delta/general}, $\delta_Q(\widetilde{C}) =
[\T^{\mathbf{m}}_n,C]_O-[\T^{\mathbf{m}}_0,C]_O-
\dim_{\mathbb{C}}(\frac{H_{\T_0}}{H_{\T_n}}) =
176-105-41=30$.
\begin{figure}
\begin{center}
\includegraphics[scale=0.6]{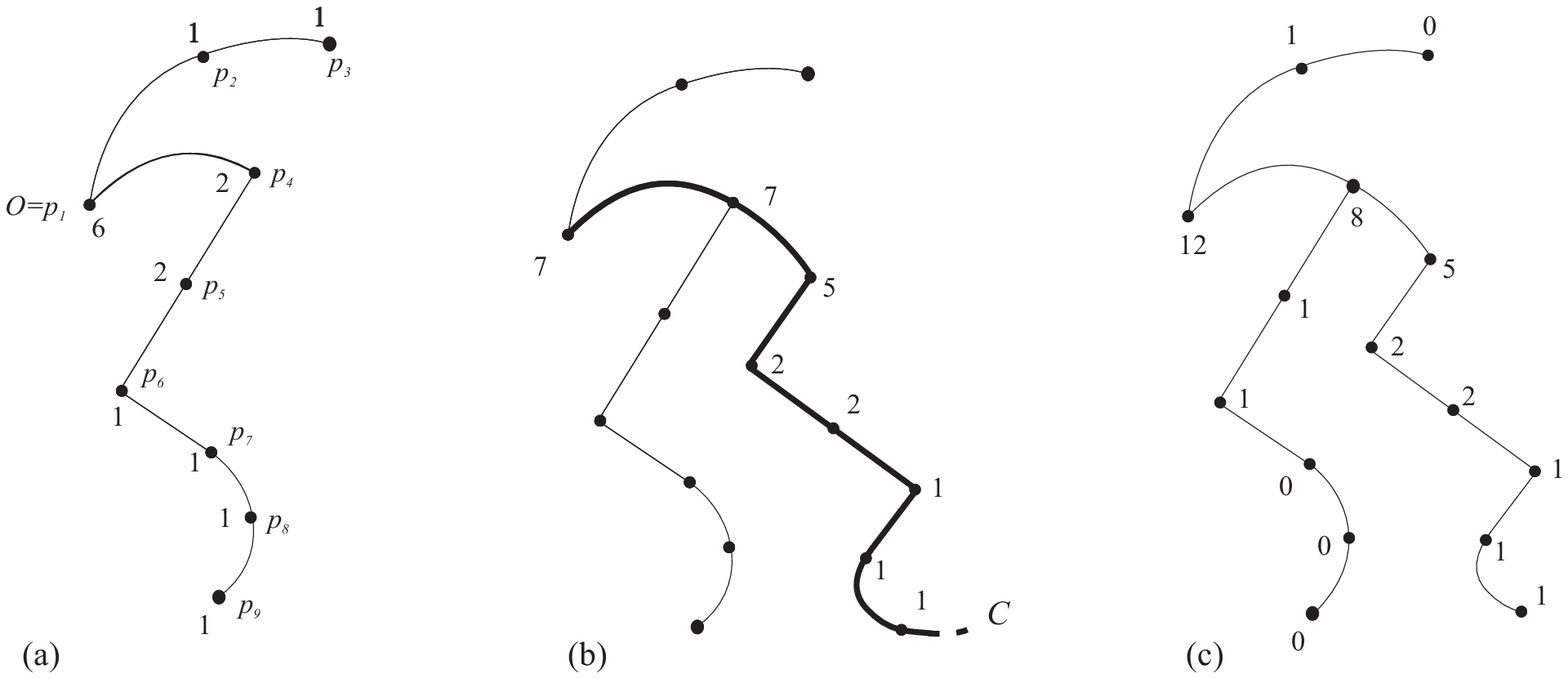}
\end{center}
\caption{\label{Fig:7}  (a) the Enriques diagram of $\K$;
(b) the singular points of $C$ are represented in bold; (c)
the Enriques diagram of $\T_n^{\boldsymbol{\omega}}$
(Example \ref{Ex:6}).}
\end{figure}
In virtue of \ref{semlemma}, we see that
\begin{eqnarray*}
\Sigma_Q(\widetilde{C}) =
\{0,7,12,14,19,21,24,26,28,31,33,35,36,38, 40,41,\ldots \\
 \ldots,42,43,45,47,48,49,50,52,53,54,55,56,57,59,60,\ldots\}
\end{eqnarray*}
Therefore, the conductor of the semigroup is $59$. Note
that the semigroup is not symmetric (see figure
\ref{Fig:7}) so, the curve $\widetilde{C}$ is not
Gorenstein (see \cite{apery, bres72}; see also
\cite{campillo}).
\end{example}

\section*{Appendix: Partial unloading}

The goal is to describe a slight modification of the
unloading procedure that we call \emph{partial unloading}.
Partial unloading has been used in section 5, but it is
explained separately for clarity and because it may be
useful in any other context where precise control of
unloading is needed. The main virtue of partial unloading
is that, fixed a non-consitent cluster $\K$ it gives rise
to an equivalent cluster, with no variation of the virtual
values at some prefixed points and with non-negative
excesses at the remaining points.

First of all, we introduce a  partial order in the set of
weighted (not necessarily consistent) clusters: if $\Q$ and
$\Q'$ are weighted clusters, we can assume that they have
the same set of points. Then, if $\Q=(Q,\tau)$ and $\Q'=(Q,
\tau')$, we write $\Q\triangleleft \Q'$ to mean that
$v_p^{\Q}\leq v_p^{\Q'}$, for each $p\in Q$. We denote by
$\mathbf{A}_K=-\textbf{P}_K^t\textbf{P}_K$ the
\emph{intersection matrix} of $E_K$, where $\textbf{P}_K$
is the \emph{proximity matrix} of $K$.

\begin{lem} \label{partunload}
Given a weighted cluster $\K=(K,\nu)$ (consistent or not)
and a subset $K_0\subset K$, there exists a cluster
$\K'=(K,\nu')$ equivalent to $\K$ and such that:
\begin{itemize}
\item[(i)] $\rho^{\K'}_p\geq 0$ if $p\in K\setminus K_0$; \item[(ii)]
$v^{\K'}_p=v^{\K}_p$ if $p\in K_0$.
\end{itemize}Moreover, with the order relation $\triangleleft$
there is a unique minimal cluster satisfying (i) and (ii).
This cluster will be denoted by $\widetilde{\K}^{K_0}$.
\end{lem}

\begin{proof}
Put $\K^{0}=\K$ and, inductively, as far as $\K^{i-1}$ has
negative excess at some point $p\in K\setminus K_0$, define
$\K^i$ from $\K^{i-1}$ by unloading on $p$. We claim that
there is an $m$ such that $\K^{m}$ has non-negative excess
at each point in $K\setminus K_0$, and
$v^{\K^m}_p=v^{\K}_p$ if $p\in K_0$. To show this, note
that the steps of the above procedure are part of an
unloading sequence giving rise to a consistent cluster as
described in Theorem 4.6.2 of \cite{CasasBook}. We reach
this cluster after finitely many steps, independently on
the order of the unloadings performed. Hence, after
finitely many steps we reach a cluster $\K^m$ satisfying
the condition (i). The condition (ii) is clear, because no
unloading is performed on any point of $K_0$.

Now, assume that $\K_{(1)}=(K,\nu^{(1)})$ and
$\K_{(2)}=(K,\nu^{(2)})$ are equivalent to $\K$ and verify
the conditions (i) and (ii). Then, for $i=1,2$ and for
every $p\in K$,
\begin{eqnarray}\label{picina}
v_p^{\K} \leq v_p^{\K^{(i)}} \leq v_p^{\widetilde{\K}}.
\end{eqnarray}
Put $\K^{(0)}=(K,\nu^{(0)})$ the cluster defined by taking
virtual values
$v^{\K^{(0)}}_p=\min\{v_p^{\K^{(1)}},v_p^{\K^{(2)}}\}$ for
every $p\in K$. Since $v_p^{\K} \leq v_p^{\K^{(0)}} \leq
v_p^{\widetilde{\K}}$, $\K^{(0)}$ is equivalent to $\K$.
For $p\in K$, write $\omega(p)=\sharp\{q\in K \mid
q\rightarrow p\}+1$. Then, if $p\in K\setminus K_0$ and say
$v_p^{\K^{(1)}}\leq v_p^{\K^{(2)}}$, we have (Artin's
trick: see Lemma 1.3 \cite{Art66})
\begin{eqnarray*}
\rho_p^{\K^{(0)}} & = &
-\textbf{1}_p^t\mathbf{A}_K\mathbf{v}^{\K^{(0)}}=
\omega(p)v^{\K^{(0)}}_p-\sum_{d(p,q)=1}v^{\K^{(0)}}_q \geq \\
& \geq &
\omega(p)v^{\K^{(1)}}_p-\sum_{d(p,q)=1}v^{\K^{(1)}}_q =
-\textbf{1}_p^t\mathbf{A}_K\mathbf{v}^{\K^{(1)}}=\rho_p^{\K^{(1)}}
\geq 0
\end{eqnarray*}
the last inequality by assumption. Moreover, if $p\in K_0$,
then $v^{\K^{(0)}}_p=v^{\K}_p$, so $\K^{(0)}$ verifies the
conditions (i) and (ii). From this, the last assertion
follows.
\end{proof}

In practice, $\widetilde{\K}^{K_0}$ is computed by
performing usual unloading on the points not in $K_0$ with
negative excess, so we say that it is  obtained from $\K$
by \emph{partial unloading (relative to $K_0$)}. Note that
if $K_0=\emptyset$, partial unloading equals usual
unloading.

\vspace{3mm} The following lemma shows that the variation
of the excesses from $\K$ to $\widetilde{\K}^{K_0}$ (also,
of the virtual multiplicities and values) is independent of
the excesses of $\K$ in the points of $K_0$. Recall that
$\mathbf{v}_{\widetilde{\K}}=\mathbf{v}_{\K}+\mathbf{n}$,
where  $\mathbf{n}$ has non-negative entries and
$\mathbf{n}=(0)$ if and only if $\K$ is consistent.

\begin{lem}\label{tech}
Let $K$ be a cluster and $K_0\subset K$ a proper subset.
Let
\[\K^i=\sum_{p\in K_0}m^i_p\K(p)+\sum_{p\in K\setminus
K_0}m_p\K(p),\quad i=1,2\] be clusters with $m_p\geq 0$ for
each $p\in K\setminus K_0$. Write $\widehat{\K^i}$ for the
cluster obtained from $\K^{i}$ by partial unloading
relative to $K_0$. Then, there exists non-negative integers
$\boldsymbol{\omega}=\{\omega_p\}_{p\in K}$ such that
$\omega_p\leq m_p$ if $p\in K\setminus K_0$ and for $i=1,2$
\[\widehat{\K^i}=\sum_{p\in
K_0}(m^i_p-\omega_p)\K(p)+\sum_{p\in K\setminus
K_0}(m_p-\omega_p)\K(p).\]
\end{lem}

\begin{proof}
By definition of partial unloading, we have
\[\textbf{v}_{\widehat{\K^i}}=\textbf{v}_{\K^i_q}+\textbf{n}^i\] where the entries of
$\textbf{n}^i=(n^i_p)_{p\in K}$ are non-negative and
$n^i_p=0$ if $p\in K_0$. Take $\textbf{n}^0=(n^0_p)_{p\in
K}$, $n^0_p=\min\{n^1_p,n^2_p\}$, and write $\K^i_0$ for
the cluster with set of points $K$ and system of values
given by $\textbf{v}_{\K^i}+\textbf{n}^0$.

I claim that $\K^i_0$ satisfies the conditions (i) and (ii)
of \ref{partunload}. If $p\in K_0$, we have that
$n^1_p=n^2_p=0$. Thus, $v_p^{\K_0^i}=v_p^{\K^i}$ and the
condition (ii) follows. For the condition (i), assume that
$p\in K\setminus K_0$ and take $j\in \{1,2\}$ such that
$n^j_p=\min\{n^1_p,n^2_p\}$. Then,
\begin{eqnarray*}
-\textbf{1}_p^t \textbf{A}_K \textbf{n}^0 = \omega(p)\
n^0_p-\sum_{d(p,q)=1}n^0_q \ \geq \ \omega(p)\
n^j_p-\sum_{d(p,q)=1}n^j_q = -\textbf{1}_p^t \textbf{A}_K
\textbf{n}^j.
\end{eqnarray*}
Hence,
\begin{eqnarray*}
\rho_p^{\K^i_0} & = & -\textbf{1}_p^t \textbf{A}_K \textbf{v}_{\K^i}-\textbf{1}_p^t \textbf{A}_K \textbf{n}^0 \geq \\
& \geq & -\textbf{1}_p^t \textbf{A}_K \textbf{v}_{\K^i}-\textbf{1}_p^t \textbf{A}_K \textbf{n}^j = \\
& = & m_p-\textbf{1}_p^t \textbf{A}_K \textbf{n}^j =
\rho_p^{\widehat{\K^{j}}}\geq 0,
\end{eqnarray*}
the last inequality by definition of $\widehat{\K}^j$. By
the minimality of $\widehat{\K^i}$, we deduce that
$\textbf{n}^0=\textbf{n}^i$ and therefore, $n^1_p=n^2_p$
for each $p\in K$. In particular, for every $p\in K$,
$\rho^{\widehat{\K^i}}_p=-\textbf{1}_p^t\textbf{A}_K(\textbf{v}_{\K^i}+\textbf{n}^0)$.
Now, if we write $\omega_p=\textbf{1}_p^t \textbf{A}_K
\textbf{n}^0$, we have that for $i=1,2$,
\begin{eqnarray*}
\rho_p^{\widehat{\K^i}} = -\textbf{1}_p^t \textbf{A}_K
\textbf{v}_{\K^i}-\omega_p=\rho_p^{\K^i}-\omega_p.
\end{eqnarray*}
Moreover, since $n^0_p=0$, we infer that
$\omega_p=\textbf{1}_p^t \textbf{A}_K
\textbf{n}^0=\sum_{d(p,q)=1}n^0_q\geq 0$.Now, if $p\in
K\setminus K_0$, $\rho_p^{\K^i}=m_p$ and so, we have
$\rho_p^{\widehat{\K^1}}=\rho_p^{\widehat{\K^2}}=m_p-\omega_p$.
Since the excess of $\widehat{\K^i}$ at $p$ is
non-negative, it is clear that $m_p\geq \omega_p$. On the
other hand, if $p\in K_0$, we have $\rho_p^{\widehat{\K^1}}
= m^1_p-\omega_p$ and $\rho_p^{\widehat{\K^2}} =
m^2_p-\omega_p$.
\end{proof}

In the study of the principality of  curves on a sandwiched
surface, partial unloding is specially useful because of
\ref{Cartier}. Keeping the notation already used , assume
that we want to prove the existence of some Cartier
divisors on $X=Bl_I(S)$ satisfying some prefixed
properties, or even, we want to construct them. The general
sketch of the procedure is the following: write
$\K=(K,\nu)$ for the cluster $BP(I)$ and take
$\K^{\mathbf{m}}=\sum_{p\in \K_+}m_p\K(p)$, where $m_p$ are
integers to be determined. Add convenient extra conditions
to $\K^{\mathbf{m}}$ (i.e. extra points with virtual
multiplicity one) in order to force the curves going
through the obtained cluster to have the desired properties
and perform partial unloading relative to $\K_+$. After
repeating this procedure a number of times, a new cluster
$\K'$ with non-negative excesses at the points of
$K\setminus \K_+$ is obtained. Moreover, for $p\in \K_+$
\begin{itemize}
\item[(i)] $\rho^{\K'}_p=m_p-\omega_p$;
\item[(ii)] $v_p^{\K'}=v_p^{\K^{\mathfrak{m}}}$.
\end{itemize}
Since $\omega_p$ is independent of $\mathbf{m}$,
$\widetilde{\K}^{\K_+}$ is consistent if we take the $m_p$
big enough. The curves going sharply through it will verify
the prefixed conditions.   The local principality of the
strict transform on $X$ of generic curves going through
this cluster follows from the condition (ii) together with
\ref{Cartier}.


\providecommand{\bysame}{\leavevmode\hbox to3em{\hrulefill}\thinspace}
\providecommand{\MR}{\relax\ifhmode\unskip\space\fi MR }
\providecommand{\MRhref}[2]{%
  \href{http://www.ams.org/mathscinet-getitem?mr=#1}{#2}
}
\providecommand{\href}[2]{#2}

\end{document}